\pdfoutput=1
\documentclass[11pt]{article}

\usepackage[margin=1in]{geometry}
\usepackage{amsmath,amssymb,amsthm}
\usepackage{mathtools}
\usepackage{bm}
\usepackage{xcolor}
\usepackage{graphicx}
\usepackage{hyperref}
\usepackage{multirow}

\hypersetup{
  colorlinks=true,
  linkcolor=blue,
  citecolor=blue,
  urlcolor=blue
}

\newtheorem{theorem}{Theorem}[section]
\newtheorem{lemma}[theorem]{Lemma}
\newtheorem{proposition}[theorem]{Proposition}
\newtheorem{corollary}[theorem]{Corollary}
\theoremstyle{definition}
\newtheorem{definition}[theorem]{Definition}
\newtheorem{assumption}[theorem]{Assumption}
\theoremstyle{remark}
\newtheorem{remark}[theorem]{Remark}

\usepackage{etoolbox}
\AtBeginEnvironment{definition}{\pushQED{\qed}}
\AtEndEnvironment{definition}{\popQED}
\AtBeginEnvironment{assumption}{\pushQED{\qed}}
\AtEndEnvironment{assumption}{\popQED}

\newcommand{\dd}{\mathrm{d}}
\newcommand{\dtalpha}{\partial_t^\alpha}
\newcommand{\dtaualpha}{\partial_\tau^\alpha}
\newcommand{\OmegaT}{\Omega_T}
\newcommand{\eps}{\varepsilon}
\newcommand{\mesh}{\mathcal{T}}
\newcommand{\energy}{E}
\newcommand{\kappaalpha}{\kappa_\alpha}

\numberwithin{equation}{section}
\title{A structure-preserving finite volume scheme for the time-fractional
Allen--Cahn equation with nondegenerate state-dependent mobility}
\author{Yonghyeon Jeon$^{1}$ \and Kiwoong Kwon$^{1}$}
\date{\small$^1$Department of Mathematical Sciences,\\
  Ulsan National Institute of Science and Technology,\\
  50 UNIST-gil, Eonyang-eup, Ulju-gun,\\
  Ulsan 44919, Republic of Korea}

\begin{document}

\maketitle

\begin{abstract}
We construct and analyze a fully discrete finite volume scheme for the
time-fractional Allen--Cahn equation with nondegenerate state-dependent
mobility. The scheme combines a standard two-point flux approximation of diffusion,
the L1 approximation of the Caputo derivative on a uniform time grid, and a
semiconvex splitting of the potential, and it discretizes the inverse mobility
by a secant average multiplying the entire L1 derivative.
The key estimate is a lower bound
for the accumulated pairing of this weighted L1 derivative with the solution
increments in terms of the accumulated L1 quadratic form. Using
this lower bound, we derive a discrete energy estimate and the
time regularity required for
compactness. For twice continuously differentiable semiconvex potentials
satisfying sign conditions outside $(-1,1)$ and initial
data in $H^1(\Omega)$
taking values in $[-1,1]$, we prove existence of at least one fully discrete
solution on every admissible finite volume mesh and every uniform time grid.
Every such solution
satisfies the maximum bound principle and a discrete energy bound relative
to the initial energy, without a time-step restriction.
In up to three space dimensions, under a mesh regularity condition with a
constant independent of the mesh size, we prove that any sequence of fully
discrete solutions has a subsequence converging strongly in
$L^2(\Omega\times(0,T))$ to a weak solution of the continuous problem as the
mesh size and time step tend to zero.
As a consequence, we obtain existence of a weak solution in
one and two space dimensions, and in three dimensions whenever such meshes
exist.
\end{abstract}

{\emergencystretch=1em\noindent\textbf{Keywords:} time-fractional Allen--Cahn equation,
finite volume method, L1 scheme, state-dependent mobility, maximum bound
principle, energy stability, convergence\par

\par}

\section{Introduction}
\label{sec:introduction}

{
{
The Allen--Cahn equation was introduced to describe antiphase boundary motion
and domain coarsening in ordered alloys \cite{1979Allen}.
}
Its unknown $u$ is a nonconserved order parameter whose values near
$\pm1$ represent the two phases, separated by diffuse transition layers of
width proportional to a parameter $\eps>0$. The equation is the gradient flow,
weighted by a nonnegative mobility $M$, of the free energy
\begin{equation*}
  \energy(u)=\int_\Omega
  \left(\frac{\eps^2}{2}|\nabla u|^2+F(u)\right)\,\dd x,
\end{equation*}
where the double-well potential $F$, for example
$F(u)=\frac14(u^2-1)^2$, favors the values $\pm1$. It reads
\begin{equation*}
  u_t=-M(u)\mu,\qquad\mu=-\eps^2\Delta u+F'(u),
\end{equation*}
where the chemical potential $\mu$ is the variational derivative of
$\energy$, and the mobility sets the rate of relaxation.
For constant mobility, after a rescaling in time, the transition layers
converge as $\eps\to0$ to interfaces moving by mean curvature
\cite{1992Evans}.
\par}

{
For sufficiently regular solutions with homogeneous Neumann boundary
conditions, this gradient flow structure gives the energy dissipation law
{
\begin{equation*}
  \frac{\dd}{\dd t}\energy(u(t))
  =-\int_\Omega M(u)|\mu|^2\,\dd x\le0.
\end{equation*}
}%
Under suitable sign conditions on $F'$, the maximum bound principle also
confines the solution to an interval containing its initial values.
These properties motivate numerical methods that preserve the maximum bound
and energy dissipation and approximate the solution accurately
\cite{2016Tang,2016Shen}.
\par}

{
For a state-dependent mobility, time discretizations that preserve the
maximum bound principle include a linear second-order backward
differentiation formula (BDF) scheme, which requires time steps bounded in terms of
the spatial mesh size and admits an $L^\infty$ error estimate \cite{2022Hou}, and a
Crank--Nicolson scheme, which under stabilization conditions that depend on
the spatial mesh also dissipates the original discrete energy without a
time-step restriction \cite{2025HouLiuJu}. Conditional a posteriori error
estimates were derived in \cite{2025BrunkGiesselmannLukacova}.
\par}

{
Time-fractional Allen--Cahn models replace the time derivative $u_t$ by a
Caputo derivative $\dtalpha$ of order $\alpha\in(0,1)$. The Caputo
derivative weights past rates of change, with greater weight given to more
recent changes, so the relaxation depends on the past evolution of the order
parameter. Numerical investigations of such models illustrate how the
fractional order affects relaxation \cite{2018Liu}. In this paper, we
consider the time-fractional Allen--Cahn equation
\begin{equation*}
  \dtalpha u=-M(u)\mu,
  \qquad \mu=-\eps^2\Delta u+F'(u),
\end{equation*}
with nondegenerate state-dependent mobility $M$.
\par}

{
The memory of the Caputo derivative affects solution regularity and the
analysis of time discretizations. Even solutions with regular initial data may have singular
time derivatives near $t=0$. For a fractional Allen--Cahn equation with
constant mobility, solution regularity and error estimates for time
discretizations were established in
\cite{2019DuYangZhou}.
\par}
 {
Temporal memory also affects the energy estimates available for time-fractional
equations. The integral dissipation estimate in \cite{2019Tang} bounds the
original energy by its initial value, but does not imply that this energy
decreases monotonically in time. Modified energies incorporating the solution
history were constructed in \cite{2022Quan,2021Fritza} and shown to decrease
in time.
\par}

{
At the discrete level, positivity results for fractional operators
\cite{2021Karaa} and discrete gradient structures on variable time grids
\cite{2024Qi} provide tools for energy estimates involving the discrete Caputo
derivative.
Related time discretizations include variable-step L1 schemes with energy
laws compatible with the classical limit $\alpha\to1$ \cite{2021Liaoa}
and higher-order linear schemes on nonuniform time grids \cite{2025Zhanga}.

For unit mobility, linear L1 schemes preserve the maximum bound and satisfy
an energy law incorporating the solution history \cite{2025Hou}. Since this
energy is the original discrete energy plus a nonnegative correction that
vanishes initially \cite[equation~(2.10)]{2025Hou}, its decay for the
predictor--corrector scheme of \cite[Theorem~2.6]{2025Hou} also bounds the
original discrete energy by its initial value. Again for unit mobility, a
discrete gradient approximation of the potential yields a scheme with a
decreasing energy incorporating the solution history \cite{2026Chen}.
For a spatial discretization of finite volume type, a scalar auxiliary
variable (SAV) finite volume element method for a time-fractional
phase-field model with constant mobility was analyzed in \cite{2026Liu},
with a bound for a modified energy containing the auxiliary variable and an
$L^2$ error estimate under regularity assumptions on the exact solution.
\par}

{
State-dependent mobility introduces an additional difficulty in the energy
estimate. In the discrete energy calculation, the inverse
mobility acts as a solution-dependent weight in the product of the L1
derivative and the current solution increment. Positivity results for a
convolution with fixed coefficients therefore do not directly establish the
nonnegativity of the sum of these weighted products over time.
\par}

{
To control these weighted products, we approximate the inverse mobility by
its average between consecutive solution values and use this coefficient to
weight the entire L1 derivative. This coefficient is the scalar mean value
discrete gradient of a primitive of the inverse mobility (see
\cite[Section~2.2]{2025Ehrhardt}). We combine this approximation with the
standard two-point flux approximation (TPFA) of diffusion, the L1 formula on
a uniform time grid, and a semiconvex splitting of the potential.
Using the classical convexity principle for matrices with nonpositive
off-diagonal entries and zero row and column sums
\cite[Proposition~1]{2004Moreau}, we derive a lower bound for the
accumulated pairing, the sum of these weighted products over time, in terms
of the accumulated L1 quadratic form.
This lower bound yields the discrete energy estimate and an estimate of the
piecewise affine time reconstructions of the discrete solutions in
$H^{(1+\alpha)/2}(0,T;L^2(\Omega))$, which gives compactness in time.
We obtain compactness in space and the limit of the diffusion term by the
finite volume arguments of \cite{2000EymardGallouetHerbin}, and we identify
the limit of the L1 derivative as a weak Caputo derivative in the sense of
\cite[Definition~2.9]{2018LiLiu}.
\par}

{
We assume that the mobility is continuously differentiable and bounded
away from zero on $[-1,1]$, and that the potential is twice continuously
differentiable, semiconvex, and satisfies {the sign conditions
\eqref{eq:potential-sign}}.
For initial data in $H^1(\Omega)$ taking values in $[-1,1]$, we prove
existence of at least one fully discrete solution on every admissible finite
volume mesh and uniform time grid. Every such solution satisfies the maximum
bound principle, an energy bound relative to the initial energy, and an
estimate for the accumulated L1 quadratic form.
Under an additional condition on the variation of the inverse mobility,
the fully discrete solution is unique.
On bounded intervals and on bounded connected polytopal
domains with Lipschitz boundary in two and
three space dimensions, under a mesh regularity condition with a constant
independent of the mesh size,
we prove that any sequence of fully discrete solutions has a subsequence
converging strongly in $L^2(\Omega\times(0,T))$ to a weak solution of
\eqref{eq:model}--\eqref{eq:initial-condition} as the spatial mesh size and
time step tend to zero.
Consequently, a weak solution of
\eqref{eq:model}--\eqref{eq:initial-condition} exists in one and two space
dimensions, and in three dimensions whenever a sequence of meshes with these
properties exists. To the best of our knowledge, this is the first existence
result for weak solutions of the time-fractional Allen--Cahn equation with a
nondegenerate state-dependent mobility. For constant mobility, the abstract
theory of time-fractional gradient flows in Hilbert spaces in
\cite{2019Akagi} applies to the time-fractional Allen--Cahn equation. For the
time-fractional Cahn--Hilliard equation with mobility, existence of weak
solutions was established in \cite{2021Fritz}.
\par}

{
{After introducing the model, the scheme, and the main results in
Section~\ref{sec:model-problem}, we prove the maximum bound principle and
the existence and conditional uniqueness of fully discrete solutions in
Section~\ref{sec:maximum-bound-principle}. In
Section~\ref{sec:energy-stability}, we prove the lower bound for the
accumulated pairing and derive the energy estimate and a priori bounds, which we
use in Section~\ref{sec:convergence-analysis} to prove the convergence
theorem. We present numerical experiments in
Section~\ref{sec:numerical_experiments} and conclude in
Section~\ref{sec:conc}. Appendix~\ref{app:technical-estimates} contains
the proofs of the auxiliary results.}
\par}

\section[Model, scheme, and main results]{{Model, scheme, and main results}}
\label{sec:model-problem}
\label{sec:method}

\subsection[Model and weak solutions]{{Model and weak solutions}}
\label{subsec:model}

{
{
{
Let $\Omega\subset\mathbb R^d$, $d\in\{1,2,3\}$, be a bounded open
interval $(a,b)$ if $d=1$ and a bounded, connected polytopal domain with
Lipschitz boundary if $d\in\{2,3\}$. We fix $T>0$, $0<\alpha<1$, and
$\eps>0$, and set $\Omega_T:=\Omega\times(0,T)$.
}
}

{
{
For an absolutely continuous function $v:[0,T]\to\mathbb R$, we define its
Caputo derivative of order $\alpha$ with initial time $0$ by
\begin{equation}
  \dtalpha v(t)
  =
  \frac{1}{\Gamma(1-\alpha)}
  \int_0^t (t-s)^{-\alpha}v'(s)\,\dd s,
  \qquad\text{for a.e. }t\in(0,T),
  \label{eq:caputo}
\end{equation}
where $\Gamma$ denotes the Euler gamma function (see
\cite[Section~1.5, Theorem~2.1(a), and equation~(2.4.17)]{2006KilbasSrivastavaTrujillo}).
}
}

{
We consider the time-fractional Allen--Cahn equation with nondegenerate
state-dependent mobility $M$,
\begin{equation}
  \dtalpha u = -M(u)\mu,
  \qquad
  \mu=-\eps^2\Delta u+F'(u)
  \qquad\text{in }\Omega_T,
  \label{eq:model}
\end{equation}
subject to the homogeneous Neumann boundary condition
\begin{equation}
  \partial_{\mathbf n}u
  :=
  \nabla u\cdot\mathbf n
  =
  0
  \qquad\text{on }\partial\Omega\times(0,T),
  \label{eq:neumann}
\end{equation}
where $\mathbf n$ denotes the outward unit normal to $\partial\Omega$, and the
initial condition
\begin{equation}
  u(\cdot,0)=u_0
  \quad\text{in }\Omega,
  \qquad
  u_0\in H^1(\Omega),
  \qquad
  |u_0|\le1
  \quad\text{a.e. in }\Omega.
  \label{eq:initial-condition}
\end{equation}
}
}

{
{
We assume that the mobility $M\in C^1([-1,1])$ is nondegenerate, that is,
there exist constants $M_*$ and $M^*$ such that
\begin{equation}
  0<M_*\le M(s)\le M^*<\infty
  \qquad\text{for every }s\in[-1,1],
  \label{eq:mobility-bounds}
\end{equation}
}%
}%
{%
{%
and that the potential $F\in C^2(\mathbb{R})$ satisfies the sign conditions
\begin{equation}
  F'(s)\le 0 \quad \text{for }s\le -1,
  \qquad
  F'(s)\ge 0 \quad \text{for }s\ge 1,
  \label{eq:potential-sign}
\end{equation}
and is semiconvex: there exists $\gamma\ge0$ such that
\begin{equation}
  F''(s)\ge -\gamma
  \qquad \text{for every }s\in\mathbb{R}.
  \label{eq:semiconvexity}
\end{equation}
With such a $\gamma$ fixed, the functions
\begin{equation}
  F_-(s)=\frac{\gamma}{2}s^2,
  \qquad
  F_+(s)=F(s)+\frac{\gamma}{2}s^2
  \label{eq:semiconvex-splitting}
\end{equation}
are convex and satisfy $F=F_+-F_-$. For example, the quartic potential
$F(s)=\frac14(s^2-1)^2$ satisfies
\eqref{eq:potential-sign}--\eqref{eq:semiconvexity} with $\gamma=1$, and the
splitting becomes $F_-(s)=\frac12s^2$ and $F_+(s)=\frac14s^4+\frac14$.
}
}

{
{
}
}

{
{
For $v\in H^1(\Omega)$ satisfying $F(v)\in L^1(\Omega)$, the free energy for
\eqref{eq:model} is
\begin{equation}
  \energy(v)
  =
  \int_\Omega
  \left(
    \frac{\eps^2}{2}|\nabla v|^2+F(v)
  \right)\,\dd x,
  \label{eq:free-energy}
\end{equation}
}%
{and for a sufficiently smooth solution with $|u|\le1$ and
the homogeneous Neumann boundary condition, the chain rule and spatial
integration by parts give, for $t>0$,}
\begin{equation}
  \frac{\dd}{\dd t}\energy(u(t))
  +
  \int_\Omega
  M(u(t))^{-1}\dtalpha u(t)\,u_t(t)\,\dd x
  =0.
  \label{eq:formal-fractional-energy-law}
\end{equation}
{Note that the second term in
\eqref{eq:formal-fractional-energy-law} has no definite sign at a fixed time,
so this identity alone does not imply that $\energy(u(t))$ is
nonincreasing.}
}

{
{
Since we do not require a weak solution to be absolutely continuous in time,
we formulate the Caputo derivative weakly by moving it onto a test function.
For $\eta\in C_c^\infty((-\infty,T);\mathbb R)$, the right-sided Caputo
derivative of order $\alpha$ with terminal time $T$ is
\begin{equation*}
  (\widetilde\partial_{T-}^\alpha\eta)(t)
  =
  -\frac1{\Gamma(1-\alpha)}
  \int_t^T(s-t)^{-\alpha}\eta'(s)\,\dd s,
  \qquad 0<t<T;
\end{equation*}
see \cite[Lemma~2.7]{2018LiLiu} for this representation. Every $v$ as in
\eqref{eq:caputo} satisfies
\begin{equation*}
  \int_0^T\dtalpha v(t)\,\eta(t)\,\dd t
  =
  \int_0^T\bigl(v(t)-v(0)\bigr)
  (\widetilde\partial_{T-}^\alpha\eta)(t)\,\dd t
\end{equation*}
by \cite[Lemma~2.8]{2018LiLiu}, and Definition~\ref{def:weak-caputo-derivative} takes the right-hand side,
which does not involve the derivative of $v$, as the definition.
}
\begin{definition}[Weak Caputo derivative]
\label{def:weak-caputo-derivative}
{
Let $u\in L^1(0,T;L^2(\Omega))$ and $u_0\in L^2(\Omega)$.
We define the weak Caputo derivative of $u$ associated with $u_0$ to be the
$L^2(\Omega)$-valued distribution $\dtalpha u$
on $(-\infty,T)$ satisfying
\begin{equation}
  \left\langle\dtalpha u,\eta\right\rangle
  =
  \int_0^T
  \bigl(u(t)-u_0\bigr)
  (\widetilde\partial_{T-}^\alpha\eta)(t)\,\dd t
  \qquad\text{in }L^2(\Omega),
  \label{eq:weak-caputo-derivative}
\end{equation}
for every $\eta\in C_c^\infty((-\infty,T);\mathbb R)$, where the integral is a
Bochner integral in $L^2(\Omega)$ (see
\cite[Definition~2.9]{2018LiLiu}). If this distribution admits a representative
in $L^2(0,T;L^2(\Omega))$, we use the same notation $\dtalpha u$ for that
representative.
}
\end{definition}
}

{
\begin{definition}[Weak solution]
\label{def:weak-solution}
{
We define a weak solution of
\eqref{eq:model}--\eqref{eq:initial-condition} to be a function
$u\in L^\infty(0,T;H^1(\Omega))\cap C([0,T];L^2(\Omega))$ such that $|u|\le1$ almost everywhere in $\OmegaT$,
$u(0)=u_0$ in $L^2(\Omega)$, the weak Caputo derivative of $u$ associated
with $u_0$ belongs to $L^2(\OmegaT)$, and
\begin{equation}
  \int_0^T\int_\Omega
  M(u)^{-1}\dtalpha u\,\varphi\,\dd x\,\dd t
  +
  \eps^2
  \int_0^T\int_\Omega
  \nabla u\cdot\nabla\varphi\,\dd x\,\dd t
  +
  \int_0^T\int_\Omega
  F'(u)\varphi\,\dd x\,\dd t
  =
  0
  \label{eq:weak-solution-formulation}
\end{equation}
holds for every test function $\varphi\in L^2(0,T;H^1(\Omega))$.
{
The regularity of $u$, the bound $|u|\le1$, the mobility bounds, and the
continuity of $F'$ on $[-1,1]$ ensure that \eqref{eq:weak-solution-formulation} is well defined.
}
The homogeneous Neumann boundary condition is imposed weakly.
}
\end{definition}
}

\begin{remark}[Inverse mobility in the weak formulation]
\label{rem:inverse-mobility-form}
The weak formulation \eqref{eq:weak-solution-formulation} tests the equation
in the form $M(u)^{-1}\dtalpha u=-\mu$ rather than $\dtalpha u=-M(u)\mu$.
For a sufficiently smooth function $u$ with $|u|\le1$ satisfying the
homogeneous Neumann boundary condition, testing the second form with $\varphi$ and integrating by
parts in space produces
\begin{equation*}
  \int_\Omega M(u)\,\Delta u\,\varphi\,\dd x
  =
  -\int_\Omega M(u)\,\nabla u\cdot\nabla\varphi\,\dd x
  -\int_\Omega M'(u)\,|\nabla u|^2\varphi\,\dd x.
\end{equation*}
Passing to the limit in the last term along a sequence of approximate
solutions would in general require strong convergence of the gradients, and
for $d\in\{2,3\}$ the term need not be defined for unbounded test functions. In
\eqref{eq:weak-solution-formulation}, the diffusion term is linear in $u$, and
the mobility appears only in the product of $M(u)^{-1}$ with the weak Caputo
derivative. Strong convergence of the approximations in $L^2(\OmegaT)$ and
weak convergence of the discrete Caputo derivatives suffice to pass to the
limit in this product (see Proposition~\ref{prop:discrete-caputo-limit}). The
scheme \eqref{eq:scheme-evolution} has the same structure, since the secant
average multiplies the L1 derivative.
\end{remark}

To describe the time regularity of the limit in
Theorem~\ref{thm:fully-discrete-convergence}, we use the following spaces.
For $0<s<1$, let $H^s(0,T;L^2(\Omega))$
be the space of all $w\in L^2(0,T;L^2(\Omega))$ with
\begin{equation*}
  [w]_{H^s(0,T;L^2(\Omega))}^2
  :=
  \int_0^T\int_0^T
  \frac{\|w(t)-w(r)\|_{L^2(\Omega)}^2}{|t-r|^{1+2s}}\,\dd r\,\dd t
  <\infty,
\end{equation*}
equipped with the norm
$\|w\|_{H^s(0,T;L^2(\Omega))}^2:=\|w\|_{L^2(0,T;L^2(\Omega))}^2
+[w]_{H^s(0,T;L^2(\Omega))}^2$.
For $0<\gamma\le1$, let $C^{0,\gamma}([0,T];L^2(\Omega))$ be the space of all
$w\colon[0,T]\to L^2(\Omega)$ with
\begin{equation*}
  \|w\|_{C^{0,\gamma}([0,T];L^2(\Omega))}
  :=
  \sup_{0\le t\le T}\|w(t)\|_{L^2(\Omega)}
  +\sup_{0\le r<t\le T}
  \frac{\|w(t)-w(r)\|_{L^2(\Omega)}}{(t-r)^\gamma}
  <\infty,
\end{equation*}
and define $H^s(0,T)$ and $C^{0,\gamma}([0,T])$ in the same way for real-valued
functions.

\subsection{Finite volume discretization}
\label{subsec:finite-volume-discretization}

{
{
Let $\mesh$ be an admissible finite volume mesh of $\Omega$. That is, $\mesh$
is a finite family of pairwise disjoint control volumes $K\subset\Omega$, each
the interior of a convex polytope, with
$\overline{\bigcup_{K\in\mesh}K}=\overline\Omega$, together with a point
$x_K\in K$ for each $K\in\mesh$, with the following properties. The faces of
$K$ are the relative interiors of the facets of $\overline K$ (the endpoints of
$K$ if $d=1$), and $\mathcal E_K$ denotes the set of faces of $K$. If $K\ne L$
and $\overline K\cap\overline L$ has positive $(d-1)$-dimensional measure,
then $\overline K\cap\overline L$ is a facet of both $\overline K$ and
$\overline L$. Its relative interior is an interior face, denoted by
$\sigma=K|L$, and the line through $x_K$ and $x_L$ is orthogonal to $\sigma$.
Every other face of $K$ is contained in $\partial\Omega$ and is called a
boundary face. For $d\in\{2,3\}$, such meshes are admissible in the sense of
\cite[Definition~10.1]{2000EymardGallouetHerbin}.
}
}

{
{
We write $\mathcal E_{\mathrm{int}}$ for the set of all interior faces and
$\mathcal E_{K,\mathrm{int}}:=\mathcal E_K\cap\mathcal E_{\mathrm{int}}$.
We write $|K|$ for the $d$-dimensional Lebesgue measure of $K$, and
$|\sigma|$ and $|\partial\Omega|$ for the $(d-1)$-dimensional Hausdorff
measures of a face $\sigma$ and of the boundary $\partial\Omega$. For
$\sigma\in\mathcal E_K$, $n_{K,\sigma}$ denotes the unit normal to $\sigma$
pointing outward from $K$. For an interior face $\sigma=K|L$, we set
$d_\sigma:=|x_K-x_L|$ and $\tau_\sigma:=|\sigma|/d_\sigma$. The numerical
normal flux is zero on each boundary face, which imposes the homogeneous
Neumann condition.
}
}

{
{
A cellwise finite volume function is a family of cell values
$v=\{v_K\}_{K\in\mesh}$, identified with the piecewise constant function
defined by $v(x)=v_K$ for almost every $x\in K$. For two such functions $v$
and $w$, the discrete $L^2$ inner product and the squared discrete $H^1$
seminorm are
\begin{equation*}
  (v,w)_{\mesh}:=\sum_{K\in\mesh}|K|v_Kw_K,
  \qquad
  |v|_{1,\mesh}^2:=
  \sum_{\sigma=K|L\in\mathcal E_{\mathrm{int}}}
  \tau_\sigma(v_L-v_K)^2,
\end{equation*}
and the two-point flux approximation (TPFA) of the Laplacian is
\begin{equation}
  (\Delta_{\mesh}v)_K
  :=
  \frac1{|K|}
  \sum_{\sigma=K|L\in\mathcal{E}_{K,\mathrm{int}}}
  \tau_\sigma (v_L-v_K),
  \qquad K\in\mesh.
  \label{eq:tpfa-laplacian}
\end{equation}
Grouping the contributions from the two cells adjacent to each interior face
$\sigma=K|L$ yields the discrete integration-by-parts identity
\begin{equation}
  -(\Delta_{\mesh}v,w)_{\mesh}
  =
  \sum_{\sigma=K|L\in\mathcal{E}_{\mathrm{int}}}
  \tau_\sigma (v_L-v_K)(w_L-w_K),
  \label{eq:discrete-ibp}
\end{equation}
and the discrete counterpart of the free energy $\energy$ is
\begin{equation}
  \energy_{\mesh}(v)
  :=
  \frac{\eps^2}{2}|v|_{1,\mesh}^2
  +
  \sum_{K\in\mesh}|K|F(v_K).
  \label{eq:discrete-free-energy}
\end{equation}
}
}

\subsection{L1 approximation of the Caputo derivative}
\label{subsec:l1-approximation}

{Given an integer $N\ge1$, we use the uniform time grid $t_n:=n\tau$,
$0\le n\le N$, with $\tau:=T/N$. For a scalar sequence
$\{v^n\}_{n=0}^N\subset\mathbb R$, let $\widehat v$ be its continuous
piecewise affine interpolant in time, that is, the continuous function on
$[0,T]$ that is affine on each interval $[t_{n-1},t_n]$ and satisfies
$\widehat v(t_n)=v^n$ for $0\le n\le N$. Since $\widehat v'$ is bounded, the
integral in \eqref{eq:caputo} for $\widehat v$ exists for every $t\in(0,T]$,
and the L1 approximation evaluates it at the grid points.
Since $\widehat v'=(v^j-v^{j-1})/\tau$ on $(t_{j-1},t_j)$, integrating the
kernel in \eqref{eq:caputo} over each interval gives, with the weights}
\begin{equation}
  a_k
  :=
  (k+1)^{1-\alpha}-k^{1-\alpha},
  \qquad k\ge0,
  \label{eq:l1-weights}
\end{equation}
{the uniform-grid L1 approximation}
\begin{equation}
  (\dtaualpha v)^n
  :=
  (\dtalpha \widehat v)(t_n)
  =
  {\kappaalpha}
  \sum_{j=1}^n a_{n-j}(v^j-v^{j-1}),
  \qquad 1\le n\le N,
  \label{eq:l1-caputo}
\end{equation}
{where}
\begin{equation}
  \kappaalpha
  :=
  \frac1{\Gamma(2-\alpha)\tau^\alpha};
  \label{eq:l1-prefactor}
\end{equation}
{see \cite[equations~(3.1)--(3.4)]{2007LinXu}.}
Since $a_k=(1-\alpha)\int_0^1(k+s)^{-\alpha}\,\dd s$ and
$s\mapsto s^{-\alpha}$ is positive, strictly decreasing, and convex on
$(0,\infty)$, the weights satisfy
\begin{equation}
  \begin{aligned}
    &a_0=1,\qquad a_k>0,\qquad a_{k+1}<a_k\quad(k\ge0),\\
    &a_{k-1}-2a_k+a_{k+1}\ge0\quad(k\ge1).
  \end{aligned}
  \label{eq:l1-weight-properties}
\end{equation}

\subsection{Fully discrete finite volume scheme}\label{subsec:method3}

{
{
{We fix an extension $q\in C^1(\mathbb R)$ of the inverse
mobility $M^{-1}$ from $[-1,1]$ to $\mathbb R$ satisfying}
\begin{equation}
  0<q_*\le q(s)\le q^*<\infty
  \qquad \text{for every }s\in\mathbb{R}.
  \label{eq:extended-inverse-mobility-bounds}
\end{equation}
{Such an extension exists because $M^{-1}$ is positive and
continuously differentiable on $[-1,1]$. For instance, if $g$ is a bounded
$C^1$ function on $\mathbb R$ that equals $-\log M$ on $[-1,1]$, then
$q:=e^{g}$ is such an extension. The primitive
$\Psi(s):=\int_0^s q(r)\,\dd r$, $s\in\mathbb R$, of $q$ satisfies $\Psi'=q$.}
}
}

{
{
{The initial value of the scheme consists of the cell averages of the
initial datum,}
\begin{equation}
  u^0:=P_{\mesh}u_0,
  \label{eq:scheme-initial-cell-average}
\end{equation}
{where $P_{\mesh}$ is the projection onto cell averages,
$(P_{\mesh}v)_K:=\frac1{|K|}\int_K v(x)\,\dd x$ for $v\in L^2(\Omega)$ and
$K\in\mesh$.}
}

{
{For $K\in\mesh$ and $1\le n\le N$, let
$\delta u_K^n:=u_K^n-u_K^{n-1}$. The secant average of the inverse mobility
is}
\begin{equation}
  \bar q_K^n
  :=
  \int_0^1
  q\!\left(u_K^{n-1}+\theta \delta u_K^n\right)
  \,\dd\theta,
  \label{eq:secant-inverse-mobility-integral}
\end{equation}
{which satisfies $q_*\le\bar q_K^n\le q^*$ by
\eqref{eq:extended-inverse-mobility-bounds}.}
The coefficient $\bar q_K^n$ is the scalar mean value discrete
gradient of $\Psi$ evaluated at $(u_K^{n-1},u_K^n)$ (see
\cite[Section~2.2]{2025Ehrhardt}).
Since $\Psi'=q$, the fundamental theorem of calculus gives the exact secant
increment identity
\begin{equation}
  \bar q_K^n \delta u_K^n
  =
  \Psi(u_K^n)-\Psi(u_K^{n-1}).
  \label{eq:secant-increment-identity}
\end{equation}
If $\delta u_K^n=0$, both sides of \eqref{eq:secant-increment-identity}
vanish, and \eqref{eq:secant-inverse-mobility-integral} gives
$\bar q_K^n=q(u_K^{n-1})$.
}
}

{
{
{
{For $1\le n\le N$ and given $u^0,\ldots,u^{n-1}$, a solution of the scheme
at step $n$ consists of cellwise finite volume functions
$u^n=\{u_K^n\}_{K\in\mesh}$ and $\mu^n=\{\mu_K^n\}_{K\in\mesh}$ such that,
for every $K\in\mesh$,}
}
\begin{align}
  \bar q_K^n(\dtaualpha u_K)^n&=-\mu_K^n,
  \label{eq:scheme-evolution}\\
  \mu_K^n
  &=
  -\eps^2(\Delta_{\mesh}u^n)_K
  +
  F_+'(u_K^n)-F_-'(u_K^{n-1}){,}
  \label{eq:scheme-chemical-potential}
\end{align}
{%
where $F_+$ and $F_-$ are the convex functions introduced through the
semiconvex splitting~\eqref{eq:semiconvex-splitting}. Evaluating $F_+'$ at
the new time level and $F_-'$ at the previous one is the convex splitting
of \cite{1998Eyre}. Since $F_+'$ is evaluated at $u_K^n$ and $\bar q_K^n$
depends on $u_K^n$ through \eqref{eq:secant-inverse-mobility-integral}, the
scheme is implicit in $u_K^n$.
}
}
}

\begin{remark}[Choice of the secant coefficient]
\label{rem:secant-coefficient-choice}
The dependence of $\bar q_K^n$ on $u_K^n$ is dictated by the energy estimate.
In Section~\ref{sec:energy-stability}, we bound
$\sum_{n=1}^m(\mu^n,\delta u^n)_{\mesh}$ for $1\le m\le N$ using
Corollary~\ref{lem:nonlinear-l1-positivity}, which applies to sums of the form
\begin{equation*}
  \sum_{n=1}^m\bigl(\Phi(v^n)-\Phi(v^{n-1})\bigr)(\dtaualpha v)^n
\end{equation*}
with one fixed function $\Phi$ satisfying $q_*\le\Phi'\le q^*$.
If a coefficient $c_K^n=c(u_K^{n-1},u_K^n)$ approximating the inverse
mobility replaced $\bar q_K^n$ in \eqref{eq:scheme-evolution}, the sum would
take this form in each cell for every sequence of cell values only if
$c_K^n\delta u_K^n=\Phi(u_K^n)-\Phi(u_K^{n-1})$ for every $n$.
Since $c_K^n$ must approximate $q(u_K^{n-1})$ as $\delta u_K^n\to0$, this
forces $\Phi'=q$, so $\Phi$ coincides with $\Psi$ up to an additive constant
and $c_K^n=\bar q_K^n$ whenever $\delta u_K^n\ne0$.
In particular, the coefficient $q(u_K^{n-1})$ satisfies
$q(u_K^{n-1})\,\delta u_K^n=\Phi(u_K^n)-\Phi(u_K^{n-1})$ for all values of
$u_K^{n-1}$ and $u_K^n$ only if $q$ is constant.
To the best of our knowledge, for the time-fractional Allen--Cahn equation
with state-dependent mobility, an energy estimate without a time-step
restriction has not been established for a scheme whose mobility coefficient
does not depend on $u_K^n$, and this question is beyond the scope of the present paper.
\end{remark}

\subsection[Main results]{{Main results}}
\label{subsec:convergence-approximation-setting}

\begin{assumption}[Spatial meshes and time grids]
\label{ass:convergence-discretizations}
Let $(\mesh_\ell)_{\ell\ge1}$ be a sequence of
admissible finite volume meshes of $\Omega$ as defined in
Section~\ref{subsec:finite-volume-discretization}. For every $\ell\ge1$,
$K\in\mesh_\ell$, and $\sigma\in\mathcal E_K$, let $d_{K,\sigma}$ denote
the distance from $x_K$ to the hyperplane containing $\sigma$. Assume that
there exists $\vartheta>0$, independent of $\ell$, such that
\begin{equation}
  d_{K,\sigma}
  \ge
  \vartheta\operatorname{diam}(K)
  \qquad
  \text{for every }\ell\ge1,\ K\in\mesh_\ell,
  \text{ and }\sigma\in\mathcal E_{K,\mathrm{int}};
  \label{eq:uniform-mesh-regularity}
\end{equation}
see \cite[Definition~10.3]{2000EymardGallouetHerbin}
for $d\in\{2,3\}$.
{For each $\ell\ge1$, we set
$h_\ell:=\max_{K\in\mesh_\ell}\operatorname{diam}(K)$, choose an integer
$N_\ell\ge1$, and define the uniform time grid by}
\begin{equation*}
  \tau_\ell=\frac{T}{N_\ell},
  \qquad
  t_\ell^n=n\tau_\ell,
  \qquad
  0\le n\le N_\ell.
\end{equation*}
{We assume that $h_\ell\to0$ and $\tau_\ell\to0$
as $\ell\to\infty$.}
\end{assumption}

For each $\ell\ge1$, let $P_\ell:=P_{\mesh_\ell}$ denote the
projection onto cell averages and set $u_\ell^0=P_\ell u_0$ as in
\eqref{eq:scheme-initial-cell-average}.
{By Theorem~\ref{thm:discrete-solvability}, we may choose a
solution of the fully discrete scheme on $(\mesh_\ell,\tau_\ell)$ and denote
it by}
\begin{equation}
  \bigl\{(u_\ell^n,\mu_\ell^n)\bigr\}_{n=1}^{N_\ell}.
  \label{eq:convergence-solution-selection}
\end{equation}

\begin{definition}[Time reconstructions]
\label{def:time-reconstructions}
{
For $1\le n\le N_\ell$, let $I_\ell^n=(t_\ell^{n-1},t_\ell^n]$. We define
the piecewise constant reconstructions $u_\ell^+$ and $u_\ell^-$ by
$u_\ell^+(t)=u_\ell^n$ and $u_\ell^-(t)=u_\ell^{n-1}$ for $t\in I_\ell^n$, and
we let $\widehat u_\ell$ be the continuous piecewise affine interpolant in
time of $u_\ell^0,\ldots,u_\ell^{N_\ell}$, taken in each cell as in
Section~\ref{subsec:l1-approximation}, so that
$\widehat u_\ell(t_\ell^n)=u_\ell^n$ for $0\le n\le N_\ell$.
}
\end{definition}

\begin{theorem}[Convergence of the fully discrete scheme]
\label{thm:fully-discrete-convergence}
Let $\Omega$, $T$, $\alpha$, and $\eps$ be as in
Section~\ref{subsec:model}. Assume that the initial datum satisfies
\eqref{eq:initial-condition}, that $M\in C^1([-1,1])$ satisfies
\eqref{eq:mobility-bounds}, and that $F\in C^2(\mathbb R)$ satisfies
\eqref{eq:potential-sign}--\eqref{eq:semiconvexity}. Let
$q\in C^1(\mathbb R)$ be {an extension of $M^{-1}$ from
$[-1,1]$ to $\mathbb R$ satisfying}
\eqref{eq:extended-inverse-mobility-bounds}, and let the sequence of spatial
meshes and time grids satisfy
Assumption~\ref{ass:convergence-discretizations}.

{For each $\ell$, we choose any fully discrete solution of
\eqref{eq:scheme-evolution}--\eqref{eq:scheme-chemical-potential} with initial
value $u_\ell^0=P_\ell u_0$, and define the reconstructions
$u_\ell^\pm$ and $\widehat u_\ell$ as in
Definition~\ref{def:time-reconstructions}.} Then there exist a subsequence
$(\ell_k)$ and a weak solution $u$ of
\eqref{eq:model}--\eqref{eq:initial-condition} such that
\begin{equation}
  u_{\ell_k}^+\longrightarrow u,\qquad
  u_{\ell_k}^-\longrightarrow u,\qquad
  \widehat u_{\ell_k}\longrightarrow u
  \quad\text{strongly in }L^2(\OmegaT).
  \label{eq:main-scheme-convergence}
\end{equation}
Moreover,
$u\in H^{(1+\alpha)/2}(0,T;L^2(\Omega))\cap C^{0,\alpha/2}([0,T];L^2(\Omega))$.
\end{theorem}

\begin{corollary}[{Existence of a weak solution}]
\label{cor:weak-solution-existence}
Let $u_0$, $M$, and $F$ satisfy the assumptions of
Theorem~\ref{thm:fully-discrete-convergence}. Then
\eqref{eq:model}--\eqref{eq:initial-condition} has at least one weak solution
if $d\in\{1,2\}$, and also if $d=3$ and there exists a sequence of admissible
finite volume meshes of $\Omega$ satisfying
Assumption~\ref{ass:convergence-discretizations}.
\end{corollary}

\begin{remark}[{Existence of mesh sequences}]
\label{rem:mesh-sequence-existence}
A sequence of meshes satisfying
Assumption~\ref{ass:convergence-discretizations} exists
for every $\Omega$ if $d\in\{1,2\}$. For $d=1$, we take uniform partitions and
let $x_K$ be the midpoint of $K$, so that \eqref{eq:uniform-mesh-regularity}
holds with $\vartheta=1/2$. For $d=2$, the domain $\Omega$ has a
triangulation into acute triangles (see \cite[Theorem~2.2]{2009Saraf}).
We subdivide every triangle repeatedly into four similar triangles by joining
the midpoints of its edges, and let $x_K$ be the circumcenter of each
triangle $K$.
These triangulations are admissible finite volume meshes (see
\cite[Example~9.1]{2000EymardGallouetHerbin}). Since their triangles are
similar to those of the initial triangulation,
\eqref{eq:uniform-mesh-regularity} holds with a constant independent of
$\ell$. For $d=3$, such a sequence exists, for example, if $\Omega$ is a box
$\prod_{i=1}^3(a_i,b_i)$: uniform Cartesian meshes with the same number of
cells in each coordinate direction, with $x_K$ the center of $K$, are
admissible and satisfy \eqref{eq:uniform-mesh-regularity} with a
constant independent of $\ell$.
\end{remark}

We prove Theorem~\ref{thm:fully-discrete-convergence} and
Corollary~\ref{cor:weak-solution-existence} in
Section~\ref{subsec:weak-formulation-convergence} after establishing the
a priori bounds and identifying the limits of the discrete terms.

\section[Maximum bound principle and discrete solvability]{{Maximum bound principle and discrete solvability}}
\label{sec:maximum-bound-principle}
\label{sec:discrete-solvability}

\subsection[Maximum bound principle]{{Maximum bound principle}}

{
{
As in the discrete maximum principle for an L1 scheme with constant mobility
in \cite[Theorem~4.1]{2019Tang} (see also \cite[Theorem~3.1]{2019DuYangZhou}
for a convolution quadrature scheme), the maximum bound principle rests on
writing the L1 approximation in terms of a convex combination of past
values. For a sequence $v^0,\ldots,v^N$ of cellwise finite volume functions,
we rewrite the L1 approximation \eqref{eq:l1-caputo}, using $a_0=1$, as
\begin{equation}
  (\dtaualpha v_K)^n
  =
  \kappaalpha\left(v_K^n-H_K^n(v)\right),
  \qquad K\in\mesh,\quad 1\le n\le N,
  \label{eq:l1-history-rewriting}
\end{equation}
where the L1 history value is
\begin{equation*}
  H_K^n(v)
  :=
  \sum_{j=1}^{n-1}
  \left(a_{n-j-1}-a_{n-j}\right)v_K^j
  +
  a_{n-1}v_K^0,
\end{equation*}
so that $H_K^1(v)=v_K^0$. Since the L1 weights are positive and strictly
decreasing by \eqref{eq:l1-weight-properties}, all coefficients in
$H_K^n(v)$ are positive, and their sum telescopes to $a_0=1$. Hence
$H_K^n(v)$ is a convex combination of $v_K^0,\ldots,v_K^{n-1}$, which gives
the convex-combination estimate
\begin{equation}
  v_K^j\in[-1,1]
  \quad\text{for }0\le j\le n-1
  \quad\Longrightarrow\quad
  H_K^n(v)\in[-1,1].
  \label{eq:l1-history-bound}
\end{equation}
}
}

{
{
For $r\in[-1,1]$, the splitting \eqref{eq:semiconvex-splitting} gives
$F_+'(s)-F_-'(r)=F'(s)+\gamma(s-r)$. By the sign condition
\eqref{eq:potential-sign}, $\gamma\ge0$, and $|r|\le1$, both terms on the
right are nonnegative for $s\ge1$ and nonpositive for $s\le-1$. Thus
\begin{align}
  F_+'(s)-F_-'(r)&\ge0&&\text{if }s\ge1,
  \label{eq:reaction-bound-upper}\\
  F_+'(s)-F_-'(r)&\le0&&\text{if }s\le-1.
  \label{eq:reaction-bound-lower}
\end{align}
}
}

{
\begin{theorem}[Maximum bound principle]
\label{thm:bound-preservation}
{
Let $\mesh$ be an admissible finite volume mesh.
{Assume that $F$ satisfies
\eqref{eq:potential-sign}--\eqref{eq:semiconvexity} and $M$ satisfies
\eqref{eq:mobility-bounds}. Let $q$ be the extension of $M^{-1}$ fixed in
\eqref{eq:extended-inverse-mobility-bounds}, and set $u^0=P_{\mesh}u_0$
as in \eqref{eq:scheme-initial-cell-average}, where $u_0$ satisfies
\eqref{eq:initial-condition}.} Then every fully discrete solution
$\bigl(\{u^n\}_{n=0}^N,\{\mu^n\}_{n=1}^N\bigr)$ of
\eqref{eq:scheme-evolution}--\eqref{eq:scheme-chemical-potential} satisfies
\begin{equation}
  u_K^n\in[-1,1]
  \qquad
  \text{for every }K\in\mesh\text{ and }0\le n\le N.
  \label{eq:bound-preservation-result}
\end{equation}
}
\end{theorem}
}

\begin{proof}
{
The cell-average initialization gives $u_K^0\in[-1,1]$.
Assume this bound holds through time level $n-1$.
By \eqref{eq:l1-history-bound}, $H_K^n(u)\in[-1,1]$, and
$\bar q_K^n\ge q_*>0$ by \eqref{eq:extended-inverse-mobility-bounds}, also
when $u_K^n\notin[-1,1]$.
{If a maximizing cell $K_+$ satisfies $u_{K_+}^n>1$, then
$(\Delta_{\mesh}u^n)_{K_+}\le0$.
Using \eqref{eq:scheme-evolution}--\eqref{eq:scheme-chemical-potential},
\eqref{eq:l1-history-rewriting}, and \eqref{eq:reaction-bound-upper}, we
obtain}
\begin{equation*}
  \begin{aligned}
    0&<\kappaalpha\bar q_{K_+}^n
        \bigl(u_{K_+}^n-H_{K_+}^n(u)\bigr)\\
     &=\eps^2(\Delta_{\mesh}u^n)_{K_+}
       -\bigl(F_+'(u_{K_+}^n)-F_-'(u_{K_+}^{n-1})\bigr)\le0,
  \end{aligned}
\end{equation*}
a contradiction.
{
If a minimizing cell $K_-$ satisfies $u_{K_-}^n<-1$, the same argument
with all inequalities reversed and \eqref{eq:reaction-bound-lower} in place
of \eqref{eq:reaction-bound-upper} gives a contradiction.
}
Hence $u_K^n\in[-1,1]$, which closes the induction.
}
\end{proof}

{
\begin{remark}[{Independence of the extension}]
\label{rem:extension-independence}
{
{
By Theorem~\ref{thm:bound-preservation}, the interval joining $u_K^{n-1}$
and $u_K^n$ lies in $[-1,1]$, on which every extension $q\in C^1(\mathbb R)$ of $M^{-1}$ satisfying
\eqref{eq:extended-inverse-mobility-bounds} equals $M^{-1}$. Hence all such extensions
give the same secant coefficients on a fully discrete solution, and these
coefficients satisfy $1/M^*\le\bar q_K^n\le1/M_*$ for $K\in\mesh$ and
$1\le n\le N$ by \eqref{eq:secant-inverse-mobility-integral} and
\eqref{eq:mobility-bounds}. Since the scheme depends on the extension only
through $\bar q_K^n$, a fully discrete solution for one such extension is a
fully discrete solution for every other, and the set of fully discrete
solutions is independent of the chosen extension.
}
}
\end{remark}
}

{
\subsection[Discrete solvability]{{Discrete solvability}}

\begin{theorem}[Discrete solvability]
\label{thm:discrete-solvability}
{
Under the assumptions of Theorem~\ref{thm:bound-preservation},
there exists at least one pair of sequences $\{u^n\}_{n=0}^N$ and
$\{\mu^n\}_{n=1}^N$ satisfying the initialization
\eqref{eq:scheme-initial-cell-average} and the fully discrete equations
\eqref{eq:scheme-evolution}--\eqref{eq:scheme-chemical-potential}.
}
\end{theorem}

\begin{proof}
{
We construct the solution by induction over the time levels. The initialization
gives $u_K^0\in[-1,1]$ for every $K\in\mesh$.
Fix $1\le n\le N$, suppose that $u^0,\ldots,u^{n-1}$ have been constructed
with $u_K^j\in[-1,1]$ for every $K\in\mesh$ and $0\le j\le n-1$, and set
$r_K:=u_K^{n-1}$. Then $r_K\in[-1,1]$, and $H_K^n(u)\in[-1,1]$ by the
convex-combination estimate \eqref{eq:l1-history-bound}. Let $J:=\#\mesh$,
identify each $v\in[-1,1]^J$ with a family of cell values
$\{v_K\}_{K\in\mesh}$, let $\bar q_K(v_K;r_K)$ be the secant average
\eqref{eq:secant-inverse-mobility-integral} with $u_K^{n-1}=r_K$ and
$u_K^n=v_K$, and define
\begin{equation*}
  \mathcal R_K^n(v)
  :=\kappaalpha\bar q_K(v_K;r_K)\bigl(v_K-H_K^n(u)\bigr)
      -\eps^2(\Delta_{\mesh}v)_K+F_+'(v_K)-F_-'(r_K).
\end{equation*}
The integral formula \eqref{eq:secant-inverse-mobility-integral} and the
continuity of $q$ and $F_+'$ show that $\mathcal R^n:[-1,1]^J\to\mathbb R^J$
is continuous.
Since $r_K,v_K\in[-1,1]$, the coefficient $\bar q_K(v_K;r_K)$ involves $q$
only on $[-1,1]$, where $q=M^{-1}$, and \eqref{eq:mobility-bounds} gives
$\bar q_K(v_K;r_K)\ge1/M^*>0$.

If $v_K=1$, then $(\Delta_{\mesh}v)_K\le0$, $v_K-H_K^n(u)\ge0$, and
\eqref{eq:reaction-bound-upper} gives $F_+'(1)-F_-'(r_K)\ge0$.
Hence $\mathcal R_K^n(v)\ge0$.
{If $v_K=-1$, then $(\Delta_{\mesh}v)_K\ge0$,
$v_K-H_K^n(u)\le0$, and \eqref{eq:reaction-bound-lower} gives
$F_+'(-1)-F_-'(r_K)\le0$. Hence $\mathcal R_K^n(v)\le0$.}
{Since $\mathcal R^n$ is continuous on $[-1,1]^J$ and, for every
$K\in\mesh$, its component $\mathcal R_K^n$ is nonnegative on the face
$\{v_K=1\}$ and nonpositive on the face $\{v_K=-1\}$, the
Poincar\'e--Miranda theorem \cite[Theorem~1]{2016Fonda} gives
$u^n\in[-1,1]^J$ with $\mathcal R^n(u^n)=0$.}

{We define $\mu^n$ by \eqref{eq:scheme-chemical-potential}.
Since $\bar q_K(u_K^n;r_K)=\bar q_K^n$ by the definition
\eqref{eq:secant-inverse-mobility-integral}, the L1 history representation
\eqref{eq:l1-history-rewriting} shows that $\mathcal R^n(u^n)=0$ is the
equation
$\bar q_K^n(\dtaualpha u_K)^n=-\mu_K^n$ in every cell.
Thus $(u^n,\mu^n)$ satisfies
\eqref{eq:scheme-evolution}--\eqref{eq:scheme-chemical-potential},
and $u^n\in[-1,1]^J$ closes the induction.}
}
\end{proof}

{
\begin{proposition}[Discrete uniqueness]
\label{prop:discrete-uniqueness}
Suppose that the assumptions of Theorem~\ref{thm:bound-preservation} hold,
set
$L_q:=\|q'\|_{L^\infty(-1,1)}$, and let $a_1=2^{1-\alpha}-1$ be the L1
weight in \eqref{eq:l1-weights}. If $a_1L_q<1/M^*$, then there exists a
unique pair of sequences $\{u^n\}_{n=0}^N$ and $\{\mu^n\}_{n=1}^N$
satisfying the initialization \eqref{eq:scheme-initial-cell-average} and
the fully discrete equations
\eqref{eq:scheme-evolution}--\eqref{eq:scheme-chemical-potential}.
\end{proposition}

\begin{proof}
Existence follows from Theorem~\ref{thm:discrete-solvability}. Let
$(u,\mu)$ and $(v,\nu)$ be two solution pairs. By
Theorem~\ref{thm:bound-preservation}, $u$ and $v$ take values in $[-1,1]$.
Fix $1\le n\le N$, assume that $u^j=v^j$ for $0\le j\le n-1$, and set
$r_K=u_K^{n-1}$ and $H_K=H_K^n(u)=H_K^n(v)$ for $K\in\mesh$.
For $n=1$, $H_K=r_K$. For $n\ge2$, the coefficient of $r_K$ in $H_K$ is
$a_0-a_1=1-a_1$, and the remaining coefficients are nonnegative with sum
$a_1$. Hence $H_K=(1-a_1)r_K+a_1\widehat H_K$ for some
$\widehat H_K\in[-1,1]$, and $|r_K-H_K|\le2a_1$ for every $n\ge1$.
For $s\in[-1,1]$, let $\bar q_K(s;r_K)$ be the secant average from the proof
of Theorem~\ref{thm:discrete-solvability}, and set
$\Phi_K(s):=\bar q_K(s;r_K)\,(s-H_K)$. The secant increment identity
\eqref{eq:secant-increment-identity} with $u_K^{n-1}=r_K$ and $u_K^n=s$
gives $\Phi_K(s)=\Psi(s)-\Psi(r_K)+(r_K-H_K)\bar q_K(s;r_K)$, and
\begin{equation*}
  \Phi_K'(s)
  =q(s)+(r_K-H_K)\int_0^1\theta\,q'\bigl(r_K+\theta(s-r_K)\bigr)\,\dd\theta
  \ge\frac1{M^*}-a_1L_q=:c_q>0,
\end{equation*}
where we used $q\ge1/M^*$ on $[-1,1]$, which follows from
\eqref{eq:mobility-bounds}. Hence
$(\Phi_K(s)-\Phi_K(t))(s-t)\ge c_q|s-t|^2$ for $s,t\in[-1,1]$.
By \eqref{eq:l1-history-rewriting} and
\eqref{eq:secant-inverse-mobility-integral}, the left-hand side of
\eqref{eq:scheme-evolution} at level $n$ equals $\kappaalpha\Phi_K(u_K^n)$
for $u$ and $\kappaalpha\Phi_K(v_K^n)$ for $v$. Set $w=u^n-v^n$.
Subtracting the two systems
\eqref{eq:scheme-evolution}--\eqref{eq:scheme-chemical-potential},
testing with $w$, and using discrete integration by parts
\eqref{eq:discrete-ibp}, we obtain
\begin{equation*}
  \kappaalpha\sum_{K\in\mesh}|K|\bigl(\Phi_K(u_K^n)-\Phi_K(v_K^n)\bigr)w_K
  +\eps^2|w|_{1,\mesh}^2
  +\sum_{K\in\mesh}|K|\bigl(F_+'(u_K^n)-F_+'(v_K^n)\bigr)w_K
  =0.
\end{equation*}
The last sum is nonnegative because $F_+$ is convex, so
$\kappaalpha c_q\|w\|_{L^2(\Omega)}^2\le0$. Hence $u^n=v^n$, and
\eqref{eq:scheme-chemical-potential} gives $\mu^n=\nu^n$. Since $u^0=v^0$
by \eqref{eq:scheme-initial-cell-average}, induction on $n$ completes the
proof.
\end{proof}

\begin{remark}[Scope of the uniqueness condition]
\label{rem:discrete-uniqueness-scope}
The condition $a_1L_q<1/M^*$ in Proposition~\ref{prop:discrete-uniqueness}
bounds $L_qM^*$, a relative variation of the inverse mobility on $[-1,1]$,
by $1/a_1$. It does not depend on the mesh, the time step, $\eps$, or the
potential, and it holds for constant mobility, since then $L_q=0$.
Since $a_1=2^{1-\alpha}-1$ decreases to $0$ as $\alpha\to1$, the
condition becomes weaker as $\alpha$ increases.
\end{remark}
}

}
{
\section[Energy stability and a priori bounds]{{Energy stability and a priori bounds}}
\label{sec:energy-stability}
}

\subsection[Nonlinear L1 positivity and energy stability]{{Nonlinear L1 positivity and energy stability}}

{
{
{
{For each integer $m\ge1$, let $L_m\in\mathbb R^{m\times m}$ be the lower
triangular matrix of the L1 weights \eqref{eq:l1-weights},}
\begin{equation}
  (L_m z)_n
  :=
  \sum_{j=1}^n a_{n-j}z_j,
  \qquad
  z\in\mathbb R^m,
  \quad 1\le n\le m,
  \label{eq:l1-convolution-matrix}
\end{equation}
and let the same formula define $L_mz\in L^2(\Omega)^m$ for
$z=(z_1,\ldots,z_m)\in L^2(\Omega)^m$.
For $v^0,\ldots,v^m$ in $\mathbb R$ or in $L^2(\Omega)$, we write
$\delta v^n:=v^n-v^{n-1}$ for $1\le n\le m$ and
$\delta v:=(\delta v^1,\ldots,\delta v^m)$. In this notation, the L1
approximation \eqref{eq:l1-caputo} of a scalar sequence is
$(\dtaualpha v)^n=\kappaalpha(L_m\delta v)_n$ for $1\le n\le m\le N$, and we
use the same formula to define $(\dtaualpha v)^n\in L^2(\Omega)$ for
$v^0,\ldots,v^m\in L^2(\Omega)$.
}
}
}
{
{
{

The following lemma is the classical convexity principle of
\cite[Proposition~1]{2004Moreau} for convex functions of the form
$\sum_{i=1}^J\Theta(v_i)$ with $\Theta:\mathbb R\to\mathbb R$ convex. We
include a direct proof for completeness.

\begin{lemma}[Convexity principle]
\label{lem:balanced-matrix-nonlinear-positivity}
Let $J\ge1$ be an integer, let $A\in\mathbb R^{J\times J}$ have nonpositive
off-diagonal entries and satisfy $A\mathbf1=A^\top\mathbf1=0$, where
$\mathbf1=(1,\ldots,1)^\top$, and let $\Phi:\mathbb R\to\mathbb R$ be
nondecreasing. Then, with $\Phi(v):=(\Phi(v_1),\ldots,\Phi(v_J))^\top$,
\begin{equation}
  \Phi(v)^\top Av\ge0
  \qquad\text{for every }v\in\mathbb R^J.
  \label{eq:balanced-matrix-nonlinear-positivity}
\end{equation}
\end{lemma}

\begin{proof}
{
Set $w_{ij}:=-A_{ij}\ge0$ for $i\ne j$ and
$\Theta(s):=\int_0^s\Phi(r)\,\dd r$.
By monotonicity of $\Phi$, we have
$\Phi(v_i)(v_i-v_j)\ge\Theta(v_i)-\Theta(v_j)$.
Therefore,
\begin{equation*}
  \begin{aligned}
    \Phi(v)^\top Av
    &=\sum_{i\ne j}w_{ij}\Phi(v_i)(v_i-v_j)\\
    &\ge\sum_{i\ne j}w_{ij}
       \bigl(\Theta(v_i)-\Theta(v_j)\bigr)=0,
  \end{aligned}
\end{equation*}
{where the first equality follows from $A\mathbf1=0$.
The final sum vanishes because
$\sum_{j\ne i}w_{ij}=A_{ii}=\sum_{j\ne i}w_{ji}$ for each $i$,
by the zero row and column sums.}
}
\end{proof}
}

\begin{lemma}[{Structure of the L1 matrix}]
\label{lem:uniform-l1-matrix-structure}
{
Let $m\ge1$ be an integer, let $D\in\mathbb R^{m\times(m+1)}$ be the forward
difference matrix, $(Dx)_p:=x_p-x_{p-1}$ for $x=(x_0,\ldots,x_m)^\top$ and
$1\le p\le m$, and let $A_m:=D^\top L_mD$, its symmetric part
$S_m:=\frac12(A_m+A_m^\top)$, and $\mathbf1:=(1,\ldots,1)^\top\in\mathbb R^{m+1}$.
Then the following properties hold.
}
\begin{enumerate}
\item[(i)] The matrix $A_m$ has nonpositive off-diagonal entries and
satisfies $A_m\mathbf1=A_m^\top\mathbf1=0$.

\item[(ii)] {For $0\le i<j\le m$, the symmetric edge weight on the set of
time indices $\{0,\ldots,m\}$,
\begin{equation}
  \omega_{ij}
  :=
  -(S_m)_{ij}
  =
  -\frac12\left((A_m)_{ij}+(A_m)_{ji}\right),
  \label{eq:l1-symmetric-edge-weight}
\end{equation}
is nonnegative. For $v^0,\ldots,v^m$ in $\mathbb R$ or in $L^2(\Omega)$,
\begin{equation}
  \sum_{n=1}^m
  \left\langle
  \delta v^n,(L_m\delta v)_n
  \right\rangle
  =
  \sum_{0\le i<j\le m}
  \omega_{ij}\|v^j-v^i\|^2,
  \label{eq:l1-quadratic-edge-representation}
\end{equation}
where $\langle\cdot,\cdot\rangle$ and $\|\cdot\|$ are the product and the
absolute value in $\mathbb R$ or the inner product and the norm of
$L^2(\Omega)$.
}

\end{enumerate}
\end{lemma}

{The proof is given in Appendix~\ref{app:l1-matrix}.}
}
}

{
{
{
By the edge representation \eqref{eq:l1-quadratic-edge-representation},
the L1 quadratic form $\sum_{n=1}^m\delta v^n(L_m\delta v)_n$ is
nonnegative. For the energy argument, we need the same inequality for scalar
sequences after the first factor $\delta v^n$ is replaced by the increment
$\Phi(v^n)-\Phi(v^{n-1})$ of a nondecreasing function
$\Phi:\mathbb R\to\mathbb R$. Applying
Lemma~\ref{lem:balanced-matrix-nonlinear-positivity} to $A_m$, we obtain
the following nonlinear L1 positivity result.
}
}
}
{
\begin{corollary}[{Nonlinear L1 positivity on a uniform grid}]
\label{lem:nonlinear-l1-positivity}
{
Let $\Phi:\mathbb R\to\mathbb R$ be nondecreasing. For every integer
$m\ge1$ and every sequence $v^0,\ldots,v^m\in\mathbb R$,
\begin{equation}
  \sum_{n=1}^m
  \left(\Phi(v^n)-\Phi(v^{n-1})\right)
  (L_m\delta v)_n
  \ge0.
  \label{eq:nonlinear-l1-positivity}
\end{equation}
Moreover, the function $\Psi$ fixed in Section~\ref{subsec:method3}
satisfies $q_*\le\Psi'\le q^*$ on $\mathbb R$ by
\eqref{eq:extended-inverse-mobility-bounds}, and
\begin{equation}
  \begin{aligned}
  q_*\sum_{n=1}^m\delta v^n(L_m\delta v)_n
  &\le
  \sum_{n=1}^m
  \left(\Psi(v^n)-\Psi(v^{n-1})\right)(L_m\delta v)_n
  \\
  &\le
  q^*\sum_{n=1}^m\delta v^n(L_m\delta v)_n.
  \end{aligned}
  \label{eq:nonlinear-l1-two-sided-bound}
\end{equation}
}
\end{corollary}
}
\begin{proof}
Let $\mathbf v=(v^0,\ldots,v^m)^\top$.
By Lemma~\ref{lem:uniform-l1-matrix-structure}(i), $A_m=D^\top L_mD$
satisfies the assumptions of
Lemma~\ref{lem:balanced-matrix-nonlinear-positivity}. Applying that lemma,
we obtain
\begin{equation}
  \begin{aligned}
  \sum_{n=1}^m
  \bigl(\Phi(v^n)-\Phi(v^{n-1})\bigr)(L_m\delta v)_n
  &=(D\Phi(\mathbf v))^\top L_mD\mathbf v\\
  &=\Phi(\mathbf v)^\top A_m\mathbf v\ge0.
  \end{aligned}
  \label{eq:nonlinear-l1-matrix-form}
\end{equation}
This proves \eqref{eq:nonlinear-l1-positivity}.
{Since $\Psi'(s)-q_*\ge0$ and $q^*-\Psi'(s)\ge0$ for every
$s\in\mathbb R$, the functions $s\mapsto\Psi(s)-q_*s$ and
$s\mapsto q^*s-\Psi(s)$ are nondecreasing.
Applying \eqref{eq:nonlinear-l1-positivity} to these functions gives the
lower and upper bounds in \eqref{eq:nonlinear-l1-two-sided-bound}, respectively.}
\end{proof}

{

{
{Corollary~\ref{lem:nonlinear-l1-positivity} controls the
accumulated pairing
$\sum_{n=1}^m\sum_{K\in\mesh}|K|\,\bar q_K^n(\dtaualpha u_K)^n\delta u_K^n$
and yields an energy bound relative to the initial energy.}
}

\begin{theorem}[{Discrete energy stability}]
\label{thm:discrete-energy-stability}
{
Under the assumptions of Theorem~\ref{thm:bound-preservation}, let
$u^0,\ldots,u^N$ be any fully discrete solution, let
$F_{\min}:=\min_{s\in[-1,1]}F(s)$, and let
$\delta u_K:=(\delta u_K^1,\ldots,\delta u_K^m)$ for $K\in\mesh$ and
$1\le m\le N$. Then, for every $1\le m\le N$,
\begin{equation}
  \begin{aligned}
  \energy_{\mesh}(u^m)
  &+
  \frac{\eps^2}{2}\sum_{n=1}^m|\delta u^n|_{1,\mesh}^2
  \\
  &+
  \kappaalpha q_*\sum_{K\in\mesh}|K|
  \sum_{n=1}^m\delta u_K^n(L_m\delta u_K)_n
  \le
  \energy_{\mesh}(u^0).
  \end{aligned}
  \label{eq:discrete-energy-dissipation}
\end{equation}
In particular,
\begin{equation}
  \energy_{\mesh}(u^m)
  \le
  \energy_{\mesh}(u^0),
  \qquad
  1\le m\le N,
  \label{eq:discrete-energy-stability}
\end{equation}
and
\begin{equation}
  \kappaalpha q_*
  \sum_{K\in\mesh}|K|
  \sum_{n=1}^m\delta u_K^n(L_m\delta u_K)_n
  \le
  \energy_{\mesh}(u^0)
  -
  |\Omega|F_{\min}.
  \label{eq:accumulated-l1-dissipation-estimate}
\end{equation}
}
\end{theorem}

\begin{proof}
{
Fix $1\le m\le N$ and $1\le n\le m$, and test
\eqref{eq:scheme-chemical-potential} with $\delta u^n$. By discrete
integration by parts \eqref{eq:discrete-ibp} and the identity
$2a(a-b)=a^2-b^2+(a-b)^2$, applied on each interior face $\sigma=K|L$ with
$a=u_L^n-u_K^n$ and $b=u_L^{n-1}-u_K^{n-1}$, we have
\begin{equation*}
  -\eps^2(\Delta_{\mesh}u^n,\delta u^n)_{\mesh}
  =\frac{\eps^2}{2}\bigl(|u^n|_{1,\mesh}^2-|u^{n-1}|_{1,\mesh}^2
    +|\delta u^n|_{1,\mesh}^2\bigr).
\end{equation*}
The convexity of $F_+$ and $F_-$ gives
$F(u_K^n)-F(u_K^{n-1})\le\bigl(F_+'(u_K^n)-F_-'(u_K^{n-1})\bigr)\delta u_K^n$
for every $K\in\mesh$. Combining these two relations with
\eqref{eq:discrete-free-energy}, we obtain
\begin{equation*}
  \energy_{\mesh}(u^n)-\energy_{\mesh}(u^{n-1})
  +\frac{\eps^2}{2}|\delta u^n|_{1,\mesh}^2
  \le(\mu^n,\delta u^n)_{\mesh}.
\end{equation*}
{Using the L1 approximation \eqref{eq:l1-caputo},
the evolution equation \eqref{eq:scheme-evolution}, and the secant increment
identity \eqref{eq:secant-increment-identity}, and applying the lower bound
in \eqref{eq:nonlinear-l1-two-sided-bound} in each cell, we obtain}
\begin{equation*}
  \begin{aligned}
    \sum_{n=1}^m(\mu^n,\delta u^n)_{\mesh}
    &=-\kappaalpha\sum_{K\in\mesh}|K|
      \sum_{n=1}^m
      \bigl(\Psi(u_K^n)-\Psi(u_K^{n-1})\bigr)
      (L_m\delta u_K)_n\\
    &\le-\kappaalpha q_*\sum_{K\in\mesh}|K|
      \sum_{n=1}^m\delta u_K^n(L_m\delta u_K)_n.
  \end{aligned}
\end{equation*}
\begin{sloppypar}
{Summing the one-step energy inequality over $n=1,\ldots,m$
and using the estimate for $\sum_{n=1}^m(\mu^n,\delta u^n)_{\mesh}$,
we obtain \eqref{eq:discrete-energy-dissipation}.
Since the L1 quadratic forms are nonnegative, 
\eqref{eq:discrete-energy-dissipation} implies \eqref{eq:discrete-energy-stability}.}
By Theorem~\ref{thm:bound-preservation}, we have
$\energy_{\mesh}(u^m)\ge|\Omega|F_{\min}$, and
\eqref{eq:accumulated-l1-dissipation-estimate} follows from
\eqref{eq:discrete-energy-dissipation}.
\end{sloppypar}
}
\end{proof}
}

\begin{remark}[{Nonuniform time grids}]
\label{rem:nonuniform-time-grids}
Theorems~\ref{thm:bound-preservation}, \ref{thm:discrete-solvability},
and~\ref{thm:discrete-energy-stability} remain valid on a time grid
$0=t_0<t_1<\cdots<t_N=T$ with variable steps if the piecewise affine
interpolant $\widehat v$ in \eqref{eq:l1-caputo} is taken on this grid. Then
\begin{equation*}
  (\dtaualpha v)^n=\sum_{j=1}^nk_{nj}(v^j-v^{j-1}),
  \qquad
  k_{nj}:=g_j(t_n),
  \qquad
  g_j(t):=\frac{1}{t_j-t_{j-1}}\int_{t_{j-1}}^{t_j}k_\alpha(t-s)\,\dd s,
\end{equation*}
and $\kappaalpha L_m$ in \eqref{eq:discrete-energy-dissipation} and
\eqref{eq:accumulated-l1-dissipation-estimate} is replaced by the lower
triangular matrix $L$ with entries $k_{nj}$, $1\le j\le n\le m$, where
$k_\alpha(r):=r^{-\alpha}/\Gamma(1-\alpha)$.
The proofs use the weights only through two properties.
First, since $k_\alpha$ is decreasing, $0<k_{n1}<k_{n2}<\cdots<k_{nn}$, so
the L1 history value $H_K^n(v)$, now defined by
$(\dtaualpha v_K)^n=k_{nn}\bigl(v_K^n-H_K^n(v)\bigr)$, is a convex combination
of $v_K^0,\ldots,v_K^{n-1}$ and satisfies \eqref{eq:l1-history-bound}. This is
the only property of the weights that
the proofs of Theorems~\ref{thm:bound-preservation}
and~\ref{thm:discrete-solvability} use.
Second, the matrix $D^\top LD$ has the properties in
Lemma~\ref{lem:uniform-l1-matrix-structure}(i). Its row and column sums
vanish because $D\mathbf1=0$. For $0\le i\le m-1$, its entry $(i,i+1)$ is
$-k_{i+1,i+1}<0$, and its entries $(i,j)$ with $j\ge i+2$ vanish.
For $i>j$, its entry $(i,j)$ equals
$k_{i+1,1}-k_{i1}$ if $j=0$ and $i<m$, $-k_{m1}$ if $(i,j)=(m,0)$,
$k_{mj}-k_{m,j+1}$ if $i=m$ and $j\ge1$, and otherwise
\begin{equation*}
  \int_{t_i}^{t_{i+1}}\bigl(g_{j+1}'(t)-g_j'(t)\bigr)\,\dd t.
\end{equation*}
The first three entries are nonpositive because $k_\alpha$ is positive and
decreasing. The
integrand is nonpositive because $-k_\alpha'$ is decreasing, so that
$-g_{j+1}'(t)\ge-k_\alpha'(t-t_j)\ge-g_j'(t)$ for $t>t_{j+1}$.
Hence Lemma~\ref{lem:balanced-matrix-nonlinear-positivity} applies to
$D^\top LD$, Corollary~\ref{lem:nonlinear-l1-positivity} holds with $L_m$
replaced by $L$, and the proof of
Theorem~\ref{thm:discrete-energy-stability} holds with $\kappaalpha L_m$
replaced by $L$. The
L1 quadratic forms are nonnegative by
\eqref{eq:balanced-matrix-nonlinear-positivity} with $\Phi(s)=s$.
The conditional uniqueness in Proposition~\ref{prop:discrete-uniqueness},
the bounds in Section~\ref{subsec:uniform-estimates-compactness},
and the convergence analysis in Section~\ref{sec:convergence-analysis} use
the explicit uniform-grid weights.
\end{remark}

\subsection{A priori bounds for the discrete solutions}
\label{subsec:uniform-estimates-compactness}

We first bound the initial discrete energies. We then combine a coercivity
estimate for the L1 quadratic form with the estimate
\eqref{eq:accumulated-l1-dissipation-estimate} for the accumulated L1
quadratic form, and test the scheme with the L1 approximation, to derive
bounds independent of the discretization.

\begin{lemma}[Estimates for the projection onto cell averages]
\label{lem:initial-projection}
Under Assumption~\ref{ass:convergence-discretizations}, for every
$\ell\ge1$ and every $v\in H^1(\Omega)$,
\begin{equation*}
  |P_\ell v|_{1,\mesh_\ell}^2
  \le
  C_{\mathrm{proj}}
  \|\nabla v\|_{L^2(\Omega)}^2,
\end{equation*}
{where $C_{\mathrm{proj}}:=\vartheta^{-2}(d/\pi^2+2/\pi)$.}
{For $u_0$ satisfying \eqref{eq:initial-condition},
the projected initial values $u_\ell^0=P_\ell u_0$ satisfy}
\begin{equation}
  |u_\ell^0|\le1
  \quad\text{a.e. in }\Omega,
  \qquad
  u_\ell^0\longrightarrow u_0
  \quad\text{strongly in }L^2(\Omega).
  \label{eq:initial-projection-properties}
\end{equation}
Moreover, their discrete energies satisfy
\begin{equation*}
  \sup_{\ell\ge1}
  \energy_{\mesh_\ell}(u_\ell^0)
  \le
  \frac{\eps^2C_{\mathrm{proj}}}{2}
  \|\nabla u_0\|_{L^2(\Omega)}^2
  +
  |\Omega|F_{\max}
  <\infty,
\end{equation*}
{where $F_{\max}:=\max_{s\in[-1,1]}F(s)$.}
\end{lemma}

{The proof is given in Appendix~\ref{app:initial-projection}.}

The scalar L1 quadratic form admits a lower bound in terms of a fractional
integral of the piecewise constant function equal to $\delta v^n$ on
$[t_{n-1},t_n)$ for $1\le n\le m$ (see
\cite[Lemma~3.1, equation~(3.9)]{2019Tang}).
{{For sequences in $L^2(\Omega)$, we use the edge
representation \eqref{eq:l1-quadratic-edge-representation} to bound the L1
quadratic form from below by a weighted sum of squared differences between
all pairs of time levels.}}

\begin{lemma}[Coercivity of the L1 quadratic form]
\label{lem:l1-fractional-coercivity}
Let $m\ge1$ be an integer, and let $v^0,\ldots,v^m$ be elements of $L^2(\Omega)$. Then
\begin{equation}
  \sum_{n=1}^m
  \langle\delta v^n,(L_m\delta v)_n\rangle_{L^2(\Omega)}
  \ge
  c_\alpha\sum_{0\le i<j\le m}
  \frac{\|v^j-v^i\|_{L^2(\Omega)}^2}{(j-i)^{2+\alpha}},
  \label{eq:l1-fractional-coercivity}
\end{equation}
where $c_\alpha:=\alpha(1-\alpha)(1+\alpha)2^{-3-\alpha}$.
\end{lemma}

{The proof is given in Appendix~\ref{app:l1-coercivity}.}

For $K\in\mesh_\ell$ and $(x,t)\in K\times I_\ell^n$, we also use the
piecewise constant reconstructions of the secant coefficient and of the
discrete Caputo derivative,
\begin{equation}
  \bar q_\ell(x,t)=\bar q_{\ell,K}^n,
  \qquad
  \mathcal D_\ell^\alpha u_\ell(x,t)=(\partial_{\tau_\ell}^\alpha u_{\ell,K})^n,
  \label{eq:discrete-caputo-reconstruction}
\end{equation}
where $u_{\ell,K}^n$ is the value of $u_\ell^n$ on $K$, and
$\bar q_{\ell,K}^n$ and $(\partial_{\tau_\ell}^\alpha u_{\ell,K})^n$ are
given by \eqref{eq:secant-inverse-mobility-integral} and \eqref{eq:l1-caputo}
on $(\mesh_\ell,\tau_\ell)$.

\begin{proposition}[A priori bounds]
\label{prop:uniform-discrete-estimates}
{
{Under the assumptions of Theorem~\ref{thm:fully-discrete-convergence},
for each $\ell\ge1$, let
$\bigl\{(u_\ell^n,\mu_\ell^n)\bigr\}_{n=1}^{N_\ell}$ be the solution
selected in \eqref{eq:convergence-solution-selection}, with the
reconstructions in Definition~\ref{def:time-reconstructions} and
\eqref{eq:discrete-caputo-reconstruction}. Let $\kappa_{\alpha,\ell}$ be the
prefactor \eqref{eq:l1-prefactor} with $\tau=\tau_\ell$, let $\delta u_\ell$
and $\delta u_{\ell,K}$, $K\in\mesh_\ell$, be the increment vectors of
$u_\ell^0,\ldots,u_\ell^{N_\ell}$ and of
$u_{\ell,K}^0,\ldots,u_{\ell,K}^{N_\ell}$ as in
Section~\ref{sec:energy-stability}, and let}
}
{
\begin{equation*}
  \begin{aligned}
  C_{\mathrm{in}}
  &:={}
  \sup_{\ell\ge1}
  \energy_{\mesh_\ell}(u_\ell^0),
  &
  C_E
  &:={}
  C_{\mathrm{in}}-|\Omega|F_{\min},
  \\
  C_\nabla
  &:={}
  \frac{2C_E}{\eps^2},
  &
  C_{\mathrm{time}}
  &:={}
  \frac{\Gamma(2-\alpha)}{c_\alpha q_*}C_E,
  \end{aligned}
\end{equation*}
{which satisfy $C_{\mathrm{in}}<\infty$ and $0\le C_E<\infty$ by
Lemma~\ref{lem:initial-projection} and
$\energy_{\mesh_\ell}(u_\ell^0)\ge|\Omega|F_{\min}$.}
}
{
{The following assertions hold for every $\ell\ge1$.}

\begin{enumerate}
\item[(i)] The cell values satisfy $|u_{\ell,K}^n|\le1$ for every
$K\in\mesh_\ell$ and $0\le n\le N_\ell$, and the discrete energy and the
discrete spatial $H^1$ seminorm satisfy
\begin{equation}
  0\le\energy_{\mesh_\ell}(u_\ell^n)-|\Omega|F_{\min}\le C_E
  \quad\text{and}\quad
  |u_\ell^n|_{1,\mesh_\ell}^2\le C_\nabla
  \qquad\text{for every }0\le n\le N_\ell.
  \label{eq:uniform-energy-gradient-estimates}
\end{equation}

\item[(ii)] The accumulated L1 quadratic form satisfies
\begin{equation}
  \kappa_{\alpha,\ell}q_*
  \sum_{K\in\mesh_\ell}|K|
  \sum_{n=1}^{N_\ell}\delta u_{\ell,K}^n(L_{N_\ell}\delta u_{\ell,K})_n
  \le C_E.
  \label{eq:uniform-accumulated-estimates}
\end{equation}

\item[(iii)] {The $L^2(\Omega)$ differences between solution values
at different time levels satisfy}
\begin{equation}
  \tau_\ell^{-\alpha}
  \sum_{0\le i<j\le N_\ell}
  \frac{
  \|u_\ell^j-u_\ell^i\|_{L^2(\Omega)}^2
  }{(j-i)^{2+\alpha}}
  \le
  C_{\mathrm{time}}{.}
  \label{eq:uniform-fractional-pair-estimate}
\end{equation}
{{In particular,}}
\begin{equation}
  \sum_{n=1}^{N_\ell}
  \|\delta u_\ell^n\|_{L^2(\Omega)}^2
  \le
  C_{\mathrm{time}}\tau_\ell^\alpha.
  \label{eq:uniform-adjacent-increment-estimate}
\end{equation}

\item[(iv)] The reconstructions $u_\ell^+$, $u_\ell^-$, and
$\widehat u_\ell$ satisfy the bounds
\begin{equation}
  \begin{aligned}
  \max_{w_\ell\in\{u_\ell^+,u_\ell^-,\widehat u_\ell\}}
  \|w_\ell\|_{L^\infty(\OmegaT)}
  &\le1,
  \\
  \max_{w_\ell\in\{u_\ell^+,u_\ell^-,\widehat u_\ell\}}
  \underset{0<t<T}{\operatorname{ess\,sup}}
  |w_\ell(t)|_{1,\mesh_\ell}^2
  &\le
  C_\nabla.
  \end{aligned}
  \label{eq:uniform-reconstruction-bounds}
\end{equation}
\item[(v)] {The discrete Caputo derivatives satisfy}
\begin{equation}
  \|\mathcal D_\ell^\alpha u_\ell\|_{L^2(\OmegaT)}^2
  \le
  C_D,
  \label{eq:discrete-caputo-uniform-l2-bound}
\end{equation}
where
$C_D:=2C_ET^{1-\alpha}/(q_*\Gamma(2-\alpha))+T|\Omega|C_F^2/q_*^2$ and
$C_F:=\max_{r,s\in[-1,1]}|F_+'(s)-F_-'(r)|$.
\end{enumerate}
}
\end{proposition}

\begin{proof}
{
Assertions (i) and (ii) follow from Lemma~\ref{lem:initial-projection}
and Theorems~\ref{thm:bound-preservation} and
\ref{thm:discrete-energy-stability}, using $F(u_{\ell,K}^n)\ge F_{\min}$
and the definitions of $C_E$ and $C_\nabla$.
Since
$\langle\delta u_\ell^n,(L_{N_\ell}\delta u_\ell)_n\rangle_{L^2(\Omega)}
=\sum_{K\in\mesh_\ell}|K|\,\delta u_{\ell,K}^n(L_{N_\ell}\delta u_{\ell,K})_n$
and $\kappa_{\alpha,\ell}^{-1}=\Gamma(2-\alpha)\tau_\ell^\alpha$, we obtain
\eqref{eq:uniform-fractional-pair-estimate} by applying
Lemma~\ref{lem:l1-fractional-coercivity} with
$v^n=u_\ell^n$ and using (ii), and \eqref{eq:uniform-adjacent-increment-estimate}
by restricting the sum to consecutive pairs.
For (iv), the bounds for $u_\ell^\pm$ follow directly from (i).
Since $\widehat u_\ell(t)$ is a convex combination of $u_\ell^{n-1}$ and
$u_\ell^n$ on each time interval, the maximum bound and convexity of the
squared discrete spatial $H^1$ seminorm give the same bounds for
$\widehat u_\ell$.
}

{For (v), fix $\ell$, omit the index $\ell$, and let
$E^j:=|u^j|_{1,\mesh}^2$. The rewriting \eqref{eq:l1-history-rewriting} is
algebraic, so for every real sequence $w^0,\ldots,w^n$ it gives
$(\dtaualpha w)^n=\kappaalpha(w^n-H^n(w))$, where $H^n(w)$ is a convex
combination of $w^0,\ldots,w^{n-1}$. The inequality
$x(x-y)\ge\frac12(x^2-y^2)$ with $x=w^n$ and $y=H^n(w)$, and Jensen's
inequality $H^n(w)^2\le H^n(w^2)$, give
\begin{equation*}
  w^n(\dtaualpha w)^n\ge\tfrac12(\dtaualpha w^2)^n.
\end{equation*}
Applying this to the differences $u_L^j-u_K^j$ across interior faces, we
obtain from \eqref{eq:discrete-ibp} and the linearity of the L1
approximation that
$-(\Delta_{\mesh}u^n,(\dtaualpha u)^n)_{\mesh}\ge\frac12(\dtaualpha E)^n$.
We multiply \eqref{eq:scheme-evolution} by $|K|(\dtaualpha u_K)^n$, sum
over $K\in\mesh$, insert \eqref{eq:scheme-chemical-potential}, and use
this inequality. Since $\bar q_K^n\ge q_*$ and $|u_K^j|\le1$ by (i),
Young's inequality then gives
\begin{equation*}
  \frac{q_*}2\|(\dtaualpha u)^n\|_{L^2(\Omega)}^2
  +\frac{\eps^2}2(\dtaualpha E)^n
  \le\frac{|\Omega|C_F^2}{2q_*},
  \qquad 1\le n\le N.
\end{equation*}
Since $\sum_{k=0}^pa_k=(p+1)^{1-\alpha}$, summation by parts gives
\begin{equation*}
  \sum_{n=1}^N\sum_{j=1}^na_{n-j}(E^j-E^{j-1})
  =E^N+\sum_{j=1}^{N-1}a_{N-j}E^j-N^{1-\alpha}E^0
  \ge-N^{1-\alpha}E^0.
\end{equation*}
We multiply the first of these two inequalities by $\tau$ and sum over $n$.
Since
$\tau\kappaalpha N^{1-\alpha}=T^{1-\alpha}/\Gamma(2-\alpha)$ and
$\eps^2E^0\le2C_E$ by (i), we obtain
\eqref{eq:discrete-caputo-uniform-l2-bound}.}
\end{proof}

\section[Convergence to a weak solution]{{Convergence to a weak solution}}
\label{sec:convergence-analysis}

{We prove Theorem~\ref{thm:fully-discrete-convergence} and
Corollary~\ref{cor:weak-solution-existence} using the a priori bounds from Section~\ref{sec:energy-stability}.
In Section~\ref{subsec:strong-convergence-discrete-solutions},
we extract a subsequence along which the three reconstructions converge
strongly to a common limit, identify its initial trace, and show that the
weak limit of the discrete spatial gradients is its spatial gradient.
In Section~\ref{subsec:passage-to-limit}, we pass to the limit in the
reaction and diffusion terms of the scheme and in the discrete Caputo
derivative weighted by the secant coefficient.}

\subsection{Strong convergence of the discrete solutions}
\label{subsec:strong-convergence-discrete-solutions}

\begin{lemma}[{$H^s$ estimate for affine interpolants}]
\label{lem:affine-fractional-time-estimate}
{
Let $N\ge1$ be an integer, and let
$\widehat v\colon[0,T]\to L^2(\Omega)$ be the continuous piecewise affine
interpolant in time of $v^0,\ldots,v^N\in L^2(\Omega)$ on the uniform time grid
$t_n=n\tau$ with $\tau=T/N$, as in Section~\ref{subsec:l1-approximation}.
Let
\begin{equation*}
  Q_\tau(v)
  =
  \tau^{-\alpha}
  \sum_{0\le i<j\le N}
  \frac{\|v^j-v^i\|_{L^2(\Omega)}^2}{(j-i)^{2+\alpha}}.
\end{equation*}
Then, for $s=(1+\alpha)/2$,
\begin{equation*}
  [\widehat v]_{H^s(0,T;L^2(\Omega))}^2
  \le
  C_{\mathrm{aff},\alpha}Q_\tau(v),
\end{equation*}
where $C_{\mathrm{aff},\alpha}>0$ depends only on $\alpha$.
}
\end{lemma}

{The proof is given in Appendix~\ref{app:affine-interpolation}.}

\begin{lemma}[{Continuous representative}]
\label{rem:fractional-time-trace}
{
Let $s=(1+\alpha)/2$. Every $w\in H^s(0,T;L^2(\Omega))$ has a continuous
representative on $[0,T]$ satisfying
\begin{equation*}
  \|w\|_{C^{0,s-1/2}([0,T];L^2(\Omega))}
  \le C_{s,T}\|w\|_{H^s(0,T;L^2(\Omega))},
\end{equation*}
where $C_{s,T}$ depends only on $s$ and $T$.
}
\end{lemma}

\begin{proof}
{
Since $s>1/2$, the scalar Sobolev embedding
$H^s(0,T)\hookrightarrow C^{0,s-1/2}([0,T])$ holds (see
\cite[Theorem~8.2, equation~(8.2)]{2012DiNezza}).
For $z$ in a countable dense subset $Z$
of the unit ball of $L^2(\Omega)$, this embedding gives a continuous
representative of the scalar function $\langle w,z\rangle_{L^2(\Omega)}$.
Since $Z$ is countable, all these representatives agree with
$\langle w,z\rangle_{L^2(\Omega)}$ outside a single null set of times.
Taking the supremum over $z\in Z$ outside this null set shows that $w$ is
H\"older continuous on a dense subset of $[0,T]$, so it extends, by the
completeness of $L^2(\Omega)$, to a continuous representative of $w$ on
$[0,T]$ satisfying the bound.
}
\end{proof}

\begin{lemma}[Spatial shift estimate]
\label{lem:finite-volume-space-change}
Let $\mesh$ be an admissible finite volume mesh of $\Omega$, let
$h:=\max_{K\in\mesh}\operatorname{diam}(K)$, and let $w$ be a cellwise finite
volume function on $\mesh$ with $|w|\le B$ for a constant $B\ge0$. Let
$\widetilde w$ be the extension of $w$ by zero to $\mathbb R^d$. Then, for
every $\xi\in\mathbb R^d$,
\begin{equation}
  \|\widetilde w(\cdot+\xi)-\widetilde w\|_{L^2(\mathbb R^d)}^2
  \le|w|_{1,\mesh}^2\,|\xi|(|\xi|+2h)+4B^2{|\partial\Omega|}\,|\xi|.
  \label{eq:zero-extension-space-change-estimate}
\end{equation}
\end{lemma}

\begin{proof}
Fix $\xi\in\mathbb R^d$. For $x,y\in\mathbb R^d$, let $[x,y]$
denote the closed line segment joining $x$ and $y$, and set
$\Omega_\xi:=\{x\in\Omega:[x,x+\xi]\subset\Omega\}$.
We first show that
\begin{equation}
  \int_{\Omega_\xi}|w(x+\xi)-w(x)|^2\,\dd x
  \le|w|_{1,\mesh}^2\,|\xi|(|\xi|+2h).
  \label{eq:space-change-interior-estimate}
\end{equation}
For $d\in\{2,3\}$, \eqref{eq:space-change-interior-estimate}
is \cite[Lemma~18.2, equation~(18.19)]{2000EymardGallouetHerbin}.
Let $d=1$, so that every face is a point and $\tau_\sigma=1/d_\sigma$.
The case $\xi=0$ is trivial, and the change of variables $y=x+\xi$ reduces
the case $\xi<0$ to the case $\xi>0$. Let $\xi>0$, and write $w_K$ for the
value of $w$ on $K$. For almost every $x\in\Omega_\xi$, the
difference $w(x+\xi)-w(x)$ is the sum of $w_L-w_K$ over the
interior faces $\sigma=K|L$ contained in $(x,x+\xi)$, where $L$ lies to the
right of $K$. By the Cauchy--Schwarz inequality,
\begin{equation*}
  |w(x+\xi)-w(x)|^2
  \le
  \Biggl(\sum_{\sigma\subset(x,x+\xi)}d_\sigma\Biggr)
  \Biggl(\sum_{\sigma=K|L\subset(x,x+\xi)}\tau_\sigma(w_L-w_K)^2\Biggr).
\end{equation*}
The first sum equals $x_{K''}-x_{K'}$, where $K'$ and $K''$
are the control volumes containing $x$ and $x+\xi$. Since
$|x_{K'}-x|\le h$ and $|x_{K''}-(x+\xi)|\le h$, it is at most
$\xi+2h$. Each interior face is contained in $(x,x+\xi)$ for $x$ in a
set of measure at most $\xi$. Integrating over $\Omega_\xi$, we obtain
\eqref{eq:space-change-interior-estimate}.

{For $d\in\{2,3\}$, the boundary construction preceding
\cite[Lemma~18.3]{2000EymardGallouetHerbin} gives
$|\Omega\setminus\Omega_{\pm\xi}|
\le|\partial\Omega|\,|\xi|$.
For $d=1$, the set $\Omega\setminus\Omega_{\pm\xi}$ is an interval of length
at most $|\xi|$, and $|\partial\Omega|=2$, so the same estimate
holds.
{Using $|a-b|^2\le2|a|^2+2|b|^2$, the zero extension,
and $|w|\le B$, we obtain}
\begin{equation*}
  \int_{\mathbb R^d\setminus\Omega_\xi}
  |\widetilde w(x+\xi)-\widetilde w(x)|^2\,\dd x
  \le2B^2\bigl(
  |\Omega\setminus\Omega_{-\xi}|+|\Omega\setminus\Omega_\xi|
  \bigr)
  \le4B^2{|\partial\Omega|}\,|\xi|.
\end{equation*}
Adding the two bounds proves the lemma.}
\end{proof}

To identify the weak limit of the discrete gradients in
Proposition~\ref{prop:strong-compactness}(iv), we use a gradient
reconstruction on pyramids. For $K\in\mesh_\ell$ and $\sigma\in\mathcal E_K$,
let $D_{K,\sigma}:=\operatorname{conv}(\{x_K\}\cup\sigma)$ be the pyramid with
apex $x_K$ and base $\sigma$ (a segment if $d=1$). For an interior face
$\sigma=K|L\in\mathcal E_{\mathrm{int}}$, let
$D_\sigma:=D_{K,\sigma}\cup D_{L,\sigma}$. The TPFA gradient
reconstruction of a cellwise finite volume function
$v_\ell=\{v_{\ell,K}\}_{K\in\mesh_\ell}$ is
\begin{equation}
  \nabla_\ell v_\ell
  =
  d\frac{v_{\ell,L}-v_{\ell,K}}{d_\sigma}n_{K,\sigma}
  \quad\text{a.e. on }D_\sigma,
  \qquad
  \sigma=K|L\in\mathcal E_{\mathrm{int}},
  \label{eq:tpfa-gradient-reconstruction}
\end{equation}
and $\nabla_\ell v_\ell:=0$ on $D_{K,\sigma}$ for every boundary face
$\sigma\in\mathcal E_K$. Reversing the orientation changes the signs of both
$v_{\ell,L}-v_{\ell,K}$ and $n_{K,\sigma}$, so $\nabla_\ell v_\ell$ is
unchanged. The pyramids $D_{K,\sigma}$, $\sigma\in\mathcal E_K$,
partition each convex control volume $K$ up to sets of measure zero, and
admissibility gives $x_L-x_K=d_\sigma n_{K,\sigma}$ and
$|D_\sigma|=|\sigma|d_\sigma/d$. Consequently,
\begin{equation}
  \sum_{\sigma\in\mathcal E_{\mathrm{int}}}|\sigma|d_\sigma
  \le d|\Omega|,
  \label{eq:interior-diamond-volume-bound}
\end{equation}
and $\|\nabla_\ell z_\ell\|_{L^2(\Omega)^d}^2=d|z_\ell|_{1,\mesh_\ell}^2$
for every cellwise finite volume function $z_\ell$.

\begin{proposition}[Compactness of the discrete solutions]
\label{prop:strong-compactness}
{
{Under the assumptions of
Theorem~\ref{thm:fully-discrete-convergence}, let $u_\ell^\pm$ and
$\widehat u_\ell$ be the reconstructions defined in
Definition~\ref{def:time-reconstructions} of the solutions selected in
\eqref{eq:convergence-solution-selection}, and let $s=(1+\alpha)/2$.} After
passing to a subsequence, which we
do not relabel, there exists a function $u$ with the following properties.
\begin{enumerate}
\item[(i)] The function $u$ has a representative satisfying
\begin{equation}
  u\in L^\infty(0,T;H^1(\Omega))
  \cap H^s(0,T;L^2(\Omega))
  \cap C^{0,\alpha/2}([0,T];L^2(\Omega)),
  \label{eq:compactness-limit-space}
\end{equation}
and
\begin{equation}
  |u|\le1\quad\text{a.e. in }\OmegaT.
  \label{eq:limit-bound-and-uniform-h1-estimate}
\end{equation}
\item[(ii)] All three reconstructions converge strongly:
\begin{equation}
  u_\ell^+\to u,\qquad u_\ell^-\to u,\qquad\widehat u_\ell\to u
  \quad\text{in }L^2(\OmegaT).
  \label{eq:strong-convergence-reconstructions}
\end{equation}
\item[(iii)] Along the same subsequence, the continuous representative
$u\colon[0,T]\to L^2(\Omega)$ satisfies
\begin{equation}
  \max_{0\le t\le T}\|\widehat u_\ell(t)-u(t)\|_{L^2(\Omega)}\to0.
  \label{eq:uniform-time-convergence}
\end{equation}
\item[(iv)] The TPFA gradients satisfy
\begin{equation}
  \nabla_\ell u_\ell^+\rightharpoonup\nabla u
  \quad\text{weakly in }L^2(\OmegaT)^d.
  \label{eq:tpfa-gradient-weak-convergence}
\end{equation}
\item[(v)] The continuous representative satisfies $u(0)=u_0$ in
$L^2(\Omega)$.
\end{enumerate}
}
\end{proposition}

\begin{proof}
{
On $I_\ell^n$, write $t=t_\ell^{n-1}+\theta\tau_\ell$ with $\theta\in(0,1]$. Then
$\widehat u_\ell(t)-u_\ell^-(t)=\theta\,\delta u_\ell^n$ and
$\widehat u_\ell(t)-u_\ell^+(t)=-(1-\theta)\,\delta u_\ell^n$. Using
$\int_0^1\theta^2\,\dd\theta=\int_0^1(1-\theta)^2\,\dd\theta=1/3$ and
\eqref{eq:uniform-adjacent-increment-estimate}, we obtain
\begin{equation}
  \max_{w_\ell\in\{u_\ell^+,u_\ell^-\}}
  \|\widehat u_\ell-w_\ell\|_{L^2(\OmegaT)}^2
  =\frac{\tau_\ell}{3}
    \sum_{n=1}^{N_\ell}\|\delta u_\ell^n\|_{L^2(\Omega)}^2
  \le\frac{C_{\mathrm{time}}}{3}\tau_\ell^{1+\alpha}.
  \label{eq:uniform-reconstruction-differences}
\end{equation}
By Lemma~\ref{lem:affine-fractional-time-estimate} and
\eqref{eq:uniform-fractional-pair-estimate},
$[\widehat u_\ell]_{H^s(0,T;L^2(\Omega))}^2\le C_{\mathrm{aff},\alpha}C_{\mathrm{time}}$,
and the maximum bound in \eqref{eq:uniform-reconstruction-bounds} gives
$\|\widehat u_\ell\|_{L^2(\OmegaT)}^2\le T|\Omega|$. Hence
$(\widehat u_\ell)$ is bounded in $H^s(0,T;L^2(\Omega))$ by a constant
independent of $\ell$.

By Lemma~\ref{rem:fractional-time-trace}, applied to the continuous
function $\widehat u_\ell$, we have, for all $t,r\in[0,T]$,
\begin{equation}
  \|\widehat u_\ell(t)-\widehat u_\ell(r)\|_{L^2(\Omega)}
  \le C_1|t-r|^{\alpha/2},
  \qquad
  C_1:=C_{s,T}\bigl(T|\Omega|+C_{\mathrm{aff},\alpha}C_{\mathrm{time}}\bigr)^{1/2}.
  \label{eq:affine-holder-estimate}
\end{equation}
Since $\widehat u_\ell(t)$ is a convex combination of $u_\ell^{n-1}$ and
$u_\ell^n$ for $t\in[t_\ell^{n-1},t_\ell^n]$,
\eqref{eq:uniform-energy-gradient-estimates} and the convexity of the squared
discrete spatial $H^1$ seminorm give $|\widehat u_\ell(t)|\le1$ and
$|\widehat u_\ell(t)|_{1,\mesh_\ell}^2\le C_\nabla$ for every $t\in[0,T]$.
Let $W_{\ell,t}$ be the extension of $\widehat u_\ell(t)$ by zero to
$\mathbb R^d$. Then $\|W_{\ell,t}\|_{L^2(\mathbb R^d)}^2\le|\Omega|$, and
Lemma~\ref{lem:finite-volume-space-change} with $B=1$ and
$h_\ell\le\operatorname{diam}(\Omega)$ gives
\begin{equation*}
  \|W_{\ell,t}(\cdot+\xi)-W_{\ell,t}\|_{L^2(\mathbb R^d)}^2
  \le C_\nabla|\xi|\bigl(|\xi|+2\operatorname{diam}(\Omega)\bigr)
     +4{|\partial\Omega|}\,|\xi|,
\end{equation*}
which tends to zero as $\xi\to0$, uniformly in $\ell$ and $t$. By the
Kolmogorov compactness lemma
\cite[Theorem~14.1]{2000EymardGallouetHerbin} with $N=d$ and
$\omega=\Omega$, the set
$\{\widehat u_\ell(t):\ell\ge1,\ 0\le t\le T\}$ is relatively compact in
$L^2(\Omega)$. Together with \eqref{eq:affine-holder-estimate}, the
Arzel\`a--Ascoli theorem gives a subsequence and a continuous
$u\colon[0,T]\to L^2(\Omega)$ satisfying
\eqref{eq:uniform-time-convergence}. In particular,
$\widehat u_\ell\to u$ in $L^2(\OmegaT)$.
{Since the right-hand side of
\eqref{eq:uniform-reconstruction-differences} tends to zero,
$u_\ell^\pm$ converge strongly to the same limit, proving
\eqref{eq:strong-convergence-reconstructions}.}
Since $\{v\in L^2(\Omega):|v|\le1\text{ a.e.}\}$ is closed in
$L^2(\Omega)$, we have $|u|\le1$. Letting $\ell\to\infty$ in
\eqref{eq:affine-holder-estimate} gives
$u\in C^{0,\alpha/2}([0,T];L^2(\Omega))$, and Fatou's lemma gives
$[u]_{H^s(0,T;L^2(\Omega))}^2\le C_{\mathrm{aff},\alpha}C_{\mathrm{time}}$.

For the spatial regularity and (iv), we use the TPFA gradient
reconstruction. For a.e.\
$t\in(0,T)$, \eqref{eq:uniform-reconstruction-bounds} and the identity
following \eqref{eq:interior-diamond-volume-bound} give
$\|\nabla_\ell u_\ell^+(t)\|_{L^2(\Omega)^d}^2\le dC_\nabla$, so the gradient
reconstructions are bounded in $L^2(\OmegaT)^d$. Let $g$ be a weak limit along a further subsequence.

Fix $\zeta\in C_c^1(\OmegaT;\mathbb R^d)$, let
$\|\nabla_x\zeta(t)\|_{L^\infty(\Omega)}$ denote the supremum over $\Omega$ of
the Frobenius norm of the Jacobian matrix of $\zeta(\cdot,t)$, and define
\begin{equation*}
  J_\ell(t)
  :=\sum_{\sigma=K|L\in\mathcal E_{\mathrm{int}}}
  (u_{\ell,L}^+(t)-u_{\ell,K}^+(t))
  \int_\sigma\zeta(y,t)\cdot n_{K,\sigma}\,\dd S_y.
\end{equation*}
Since $\zeta$ vanishes near $\partial\Omega$, grouping the contributions of
adjacent cells and applying the divergence theorem on each cell gives
\begin{equation*}
  J_\ell(t)
  =-\int_\Omega u_\ell^+(x,t)\operatorname{div}_x\zeta(x,t)\,\dd x.
\end{equation*}
For $x\in D_{K,\sigma}$ and $y\in\sigma$, both points belong to
$\overline K$, so $|x-y|\le h_\ell$. Averaging over $x\in D_\sigma$ and
$y\in\sigma$ gives
\begin{equation*}
  \left|
  \frac1{|D_\sigma|}\int_{D_\sigma}\zeta(x,t)\,\dd x
  -\frac1{|\sigma|}\int_\sigma\zeta(y,t)\,\dd S_y
  \right|
  \le h_\ell\|\nabla_x\zeta(t)\|_{L^\infty(\Omega)}.
\end{equation*}
By the Cauchy--Schwarz inequality and
\eqref{eq:interior-diamond-volume-bound},
\begin{equation*}
  \left|
  \int_\Omega\nabla_\ell u_\ell^+(t)\cdot\zeta(t)\,\dd x-J_\ell(t)
  \right|
  \le h_\ell\sqrt{d|\Omega|}\,
  |u_\ell^+(t)|_{1,\mesh_\ell}
  \|\nabla_x\zeta(t)\|_{L^\infty(\Omega)}.
\end{equation*}
{Integrating in time and passing to the limit with $h_\ell\to0$,
the bound \eqref{eq:uniform-reconstruction-bounds}, the weak convergence
$\nabla_\ell u_\ell^+\rightharpoonup g$ in $L^2(\OmegaT)^d$, and the
strong convergence $u_\ell^+\to u$ in $L^2(\OmegaT)$, we obtain}
\begin{equation*}
  \int_{\OmegaT}g\cdot\zeta\,\dd x\,\dd t
  =-\int_{\OmegaT}u\operatorname{div}_x\zeta\,\dd x\,\dd t.
\end{equation*}
{%
Thus $u\in L^2(0,T;H^1(\Omega))$ and $g=\nabla u$, which proves
\eqref{eq:tpfa-gradient-weak-convergence}.

The bound on $\nabla_\ell u_\ell^+(t)$ passes to $\nabla u(t)$ by weak
lower semicontinuity on measurable time subsets, so
$\nabla u\in L^\infty(0,T;L^2(\Omega)^d)$. Since $|u|\le1$, this gives
$u\in L^\infty(0,T;H^1(\Omega))$.
}%
Finally, \eqref{eq:uniform-time-convergence} gives
$\widehat u_\ell(0)\to u(0)$ in $L^2(\Omega)$. Since
$\widehat u_\ell(0)=u_\ell^0\to u_0$ in $L^2(\Omega)$ by
Lemma~\ref{lem:initial-projection}, this proves (v).
}
\end{proof}

\subsection{Passage to the limit in the scheme}
\label{subsec:passage-to-limit}

Throughout this subsection, we fix the subsequence and the limit $u$ given
by Proposition~\ref{prop:strong-compactness}, and all limits as $\ell\to\infty$
are taken along this subsequence.

\subsubsection{The reaction and diffusion terms}
\label{subsec:limits-nonlinear-diffusion}

The reaction term $F_+'(u_K^n)-F_-'(u_K^{n-1})$ in
\eqref{eq:scheme-chemical-potential} has the piecewise constant
reconstruction $G_\ell:=F_+'(u_\ell^+)-F_-'(u_\ell^-)$.
{%
By Remark~\ref{rem:extension-independence}, the interval joining
$u_{\ell,K}^{n-1}$ and $u_{\ell,K}^n$ lies in $[-1,1]$, where $q=M^{-1}$, for
every $K\in\mesh_\ell$ and $1\le n\le N_\ell$, and $|u|\le1$ almost
everywhere in $\OmegaT$ by \eqref{eq:limit-bound-and-uniform-h1-estimate}.
Since $F_\pm'$ and $q$ are Lipschitz on $[-1,1]$, we have
\begin{equation}
  |G_\ell-F'(u)|+|\bar q_\ell-M(u)^{-1}|
  \le C\bigl(|u_\ell^+-u|+|u_\ell^--u|\bigr),
  \label{eq:reaction-secant-lipschitz}
\end{equation}
where $C$ is independent of $\ell$.
{By Proposition~\ref{prop:strong-compactness}, the right-hand
side of \eqref{eq:reaction-secant-lipschitz} tends to zero in
$L^2(\OmegaT)$.
Therefore,}
}
\begin{align}
  G_\ell
  &\longrightarrow
  F'(u)
  \quad\text{strongly in }L^2(\OmegaT),
  \label{eq:strong-reaction-convergence}
  \\
  \bar q_\ell
  &\longrightarrow
  M(u)^{-1}
  \quad\text{strongly in }L^2(\OmegaT).
  \label{eq:strong-secant-coefficient-convergence}
\end{align}

We test the scheme with the time averages of
$\varphi\in C_c^\infty(\mathbb R^d\times(0,T))$ at the points $x_K$,
\begin{equation}
  \varphi_{\ell,K}^n
  :=\frac1{\tau_\ell}\int_{I_\ell^n}\varphi(x_K,t)\,\dd t,
  \label{eq:test-center-time-average}
\end{equation}
and let $\varphi_\ell$ be the piecewise constant function equal to
$\varphi_{\ell,K}^n$ on $K\times I_\ell^n$.
Since $|x-x_K|\le h_\ell$ and $|t-s|\le\tau_\ell$ for $t,s\in I_\ell^n$,
we obtain
\begin{equation}
  \|\varphi_\ell-\varphi\|_{L^\infty(\OmegaT)}
  \le h_\ell\|\nabla_x\varphi\|_{L^\infty(\mathbb R^d\times(0,T))}
     +\tau_\ell\|\partial_t\varphi\|_{L^\infty(\mathbb R^d\times(0,T))}
  \longrightarrow0.
  \label{eq:test-reconstruction-uniform-convergence}
\end{equation}

\begin{proposition}[Limit of the discrete diffusion term]
\label{prop:tpfa-diffusion-limit}
{
For every $\varphi\in C_c^\infty(\mathbb R^d\times(0,T))$, with
the test values $\varphi_{\ell,K}^n$ defined by
\eqref{eq:test-center-time-average}, we have
\begin{equation}
  \begin{aligned}
  &\eps^2\tau_\ell
  \sum_{n=1}^{N_\ell}
  \sum_{\sigma=K|L\in\mathcal E_{\mathrm{int}}}
  \tau_\sigma
  \left(u_{\ell,L}^n-u_{\ell,K}^n\right)
  \left(\varphi_{\ell,L}^n-\varphi_{\ell,K}^n\right)
  \\
  &\qquad\longrightarrow
  \eps^2
  \int_0^T\int_\Omega
  \nabla u\cdot\nabla\varphi\,\dd x\,\dd t.
  \end{aligned}
  \label{eq:tpfa-diffusion-limit}
\end{equation}
}
\end{proposition}

\begin{proof}
{
Fix $t\in(0,T)$ and a cellwise finite volume function $v_\ell$. Let
$\|D_x^2\varphi(t)\|_{L^\infty(\mathbb R^d)}$ denote the supremum over
$\mathbb R^d$ of the Frobenius norm of the Hessian matrix of
$\varphi(\cdot,t)$, and let $\|D_x^2\varphi\|_{L^\infty(\mathbb R^d\times(0,T))}$
denote its supremum over $t\in(0,T)$.
Since $x_L-x_K=d_\sigma n_{K,\sigma}$,
\begin{equation*}
  \frac{\varphi(x_L,t)-\varphi(x_K,t)}{d_\sigma}
  =\int_0^1\nabla\varphi(x_K+r(x_L-x_K),t)
  \cdot n_{K,\sigma}\,\dd r.
\end{equation*}
Since every point of $[x_K,x_L]$ is within $2h_\ell$ of every point of
$D_\sigma$, we compare the averages of $\nabla\varphi(\cdot,t)\cdot n_{K,\sigma}$
over $[x_K,x_L]$ and $D_\sigma$, and the Cauchy--Schwarz inequality and
\eqref{eq:interior-diamond-volume-bound} give
\begin{equation}
  \begin{aligned}
  &\left|
  \sum_{\sigma=K|L\in\mathcal E_{\mathrm{int}}}
  \tau_\sigma(v_{\ell,L}-v_{\ell,K})
  (\varphi(x_L,t)-\varphi(x_K,t))
  -\int_\Omega\nabla_\ell v_\ell\cdot\nabla\varphi(t)\,\dd x
  \right|
  \\
  &\qquad\le
  2h_\ell\sqrt{d|\Omega|}\,|v_\ell|_{1,\mesh_\ell}
  \|D_x^2\varphi(t)\|_{L^\infty(\mathbb R^d)};
  \end{aligned}
  \label{eq:tpfa-flux-consistency}
\end{equation}
for a related flux consistency argument, see
\cite[proof of Theorem~10.3, equations~(10.36)--(10.37)]{2000EymardGallouetHerbin}.

By \eqref{eq:test-center-time-average}, the left-hand side of
\eqref{eq:tpfa-diffusion-limit} is $\eps^2$ times the time integral of the
sum over interior faces in \eqref{eq:tpfa-flux-consistency} with
$v_\ell=u_\ell^+(t)$. Integrating \eqref{eq:tpfa-flux-consistency} in time
and using \eqref{eq:uniform-reconstruction-bounds}, we bound the difference
between this integral and
$\int_{\OmegaT}\nabla_\ell u_\ell^+\cdot\nabla\varphi\,\dd x\,\dd t$ by
$2h_\ell T\sqrt{d|\Omega|C_\nabla}\,
\|D_x^2\varphi\|_{L^\infty(\mathbb R^d\times(0,T))}$, which tends to zero as
$h_\ell\to0$. Using \eqref{eq:tpfa-gradient-weak-convergence} and
multiplying by $\eps^2$, we obtain \eqref{eq:tpfa-diffusion-limit}.
}
\end{proof}

\subsubsection{The discrete Caputo derivative}
\label{subsec:discrete-caputo-identification}

\begin{proposition}[{{Convergence of the discrete Caputo derivative}}]
\label{prop:discrete-caputo-limit}
{
The following assertions hold.

\begin{enumerate}
\item[(i)] We have
\begin{equation}
  \mathcal D_\ell^\alpha u_\ell
  \rightharpoonup
  \dtalpha u
  \quad\text{weakly in }L^2(\OmegaT),
  \label{eq:discrete-caputo-weak-limit}
\end{equation}
where $\dtalpha u$ is the weak Caputo derivative of $u$ associated with
$u_0$. In particular,
$\dtalpha u\in L^2(\OmegaT)$.

\item[(ii)] For every $\varphi\in C_c^\infty(\mathbb R^d\times(0,T))$, with
$\varphi_\ell$ defined after \eqref{eq:test-center-time-average},
\begin{equation}
  \int_{\OmegaT}
  \bar q_\ell\mathcal D_\ell^\alpha u_\ell\varphi_\ell
  \,\dd x\,\dd t
  \longrightarrow
  \int_{\OmegaT}
  M(u)^{-1}\dtalpha u\,\varphi
  \,\dd x\,\dd t.
  \label{eq:weighted-discrete-caputo-limit}
\end{equation}
\end{enumerate}
}
\end{proposition}

\begin{proof}
{
Since $u_\ell^+\to u$ in $L^2(\OmegaT)$ by
Proposition~\ref{prop:strong-compactness}, $u_\ell^0\to u_0$ in
$L^2(\Omega)$ by Lemma~\ref{lem:initial-projection}, and
\eqref{eq:discrete-caputo-uniform-l2-bound} holds,
Lemma~\ref{lem:hilbert-l1-limit-identification} shows that
$\mathcal D_\ell^\alpha u_\ell$ converges weakly in $L^2(\OmegaT)$ to the
representative of the weak Caputo derivative of $u$ associated with $u_0$.
This proves (i). For (ii), we write
$\bar q_\ell\varphi_\ell-M(u)^{-1}\varphi
=\bar q_\ell(\varphi_\ell-\varphi)+(\bar q_\ell-M(u)^{-1})\varphi$.
Since $\bar q_\ell\le q^*$ and $\varphi$ is bounded,
\eqref{eq:strong-secant-coefficient-convergence} and
\eqref{eq:test-reconstruction-uniform-convergence} give
$\bar q_\ell\varphi_\ell\to M(u)^{-1}\varphi$ strongly in $L^2(\OmegaT)$.
Pairing this with \eqref{eq:discrete-caputo-weak-limit} gives
\eqref{eq:weighted-discrete-caputo-limit}.
}
\end{proof}

\subsubsection{Proofs of the main results}
\label{subsec:weak-formulation-convergence}

\begin{proof}[Proof of Theorem~\ref{thm:fully-discrete-convergence}]
{
{
Proposition~\ref{prop:strong-compactness} gives a subsequence and a function
$u$ with the regularity in \eqref{eq:compactness-limit-space}, all three strong
convergences in \eqref{eq:main-scheme-convergence}, $|u|\le1$, and
$u(0)=u_0$ in $L^2(\Omega)$.
Along this subsequence, Proposition~\ref{prop:discrete-caputo-limit}(i) shows
that the weak Caputo derivative of $u$ associated with $u_0$ belongs to
$L^2(\OmegaT)$.
}

{
To verify \eqref{eq:weak-solution-formulation}, we fix
$\varphi\in C_c^\infty(\mathbb R^d\times(0,T))$ and use the test values
\eqref{eq:test-center-time-average}.
} We eliminate $\mu_\ell^n$
from the scheme, multiply by $\tau_\ell|K|\varphi_{\ell,K}^n$, and sum over
cells and time levels. By discrete integration by parts, we obtain
\begin{equation*}
  \begin{aligned}
  0={}&\int_{\OmegaT}\bar q_\ell\mathcal D_\ell^\alpha u_\ell
             \varphi_\ell\,\dd x\,\dd t\\
     &+\eps^2\tau_\ell\sum_{n=1}^{N_\ell}
        \sum_{\sigma=K|L\in\mathcal E_{\mathrm{int}}}
        \tau_\sigma(u_{\ell,L}^n-u_{\ell,K}^n)
                    (\varphi_{\ell,L}^n-\varphi_{\ell,K}^n)\\
     &+\int_{\OmegaT}G_\ell\varphi_\ell\,\dd x\,\dd t.
  \end{aligned}
\end{equation*}
{%
Passing to the limit in the three terms on the right-hand side by
\eqref{eq:weighted-discrete-caputo-limit}, by \eqref{eq:tpfa-diffusion-limit},
and by \eqref{eq:strong-reaction-convergence} and
\eqref{eq:test-reconstruction-uniform-convergence}, respectively, we obtain
}
\begin{equation}
  \int_{\OmegaT}
  \left(M(u)^{-1}\dtalpha u\,\varphi
       +\eps^2\nabla u\cdot\nabla\varphi
       +F'(u)\varphi\right)\,\dd x\,\dd t=0.
  \label{eq:limit-weak-formulation}
\end{equation}

{
{Since $\Omega$ has Lipschitz boundary, the restrictions to
$\Omega$ of functions in $C_c^\infty(\mathbb R^d)$ are dense in
$H^1(\Omega)$ (see \cite[Theorem~3.22]{2003Adams}).
Finite sums of products of these spatial functions and functions in
$C_c^\infty(0,T)$ are therefore dense in $L^2(0,T;H^1(\Omega))$.}
The established $L^2$ bounds for $M(u)^{-1}\dtalpha u$,
$\nabla u$, and $F'(u)$ show that the left-hand side of
\eqref{eq:limit-weak-formulation} is a continuous linear functional on
$L^2(0,T;H^1(\Omega))$, so the identity holds for every test function
in this space by density.
} The regularity \eqref{eq:compactness-limit-space}, the initial value
$u(0)=u_0$, the bound $|u|\le1$, and $\dtalpha u\in L^2(\OmegaT)$ show that
$u$ satisfies Definition~\ref{def:weak-solution}.
}
\end{proof}

\begin{proof}[Proof of Corollary~\ref{cor:weak-solution-existence}]
By Remark~\ref{rem:mesh-sequence-existence} if $d\in\{1,2\}$, and by
assumption if $d=3$, there is a sequence of meshes and time grids satisfying
Assumption~\ref{ass:convergence-discretizations}. We fix an extension $q$ of
$M^{-1}$ as in \eqref{eq:extended-inverse-mobility-bounds}. Then
Theorem~\ref{thm:discrete-solvability} gives, for each $\ell$, a fully
discrete solution with initial value $P_\ell u_0$, and
Theorem~\ref{thm:fully-discrete-convergence} gives a subsequence converging
to a weak solution.
\end{proof}

{
\section[Numerical experiments]{{Numerical experiments}}
\label{sec:numerical_experiments}
%=========================================================================

In this section, we present numerical experiments in two and three dimensions that examine the accuracy and the structure-preserving properties of the fully discrete finite volume scheme.

We consider three mobilities satisfying the nondegeneracy assumption \eqref{eq:mobility-bounds}: a constant mobility $M_0$, a weakly varying mobility $M_1$, and a regularized degenerate mobility $M_2$ that varies strongly on $[-1,1]$.
Table~\ref{tab:mobility_models} lists them together with the inverse mobility $q=M^{-1}$ and its primitive $\Psi$ from Section~\ref{subsec:method3} on $[-1,1]$.
Using these explicit formulas and \eqref{eq:Psi_2}, we evaluate the secant average \eqref{eq:secant-inverse-mobility-integral} through the identity \eqref{eq:secant-increment-identity} without numerical quadrature.
\begin{table}[htbp]
    \centering
    \begin{tabular}{lccc}
        \hline
        Model & $M(u)$ & $q(u)$ & $\Psi(u)$ \\
        \hline
        $M_0$: constant & $1$ & $1$ & $u$ \\
        $M_1$: weakly varying & $1 - 0.5u^2$ & $\frac{2}{2-u^2}$ & $\frac{1}{\sqrt{2}} \ln\left| \frac{\sqrt{2}+u}{\sqrt{2}-u} \right|$ \\
        $M_2$: regularized degenerate & $(1 - u^2)^2 + 10^{-10}$ & $\frac{1}{(1 - u^2)^2 + 10^{-10}}$ & see \eqref{eq:Psi_2} \\
        \hline
    \end{tabular}
    \caption{Mobility models, their inverses $q$ on $[-1,1]$, and the primitives $\Psi$. The regularization constant $10^{-10}$ gives $M_2\ge10^{-10}>0$.}
    \label{tab:mobility_models}
\end{table}

For $M_2$, the primitive is
\begin{equation} \label{eq:Psi_2}
    \Psi_2(u) = \frac{1}{4AB} \ln\left(\frac{u^2+Bu+A}{u^2-Bu+A}\right) + \frac{1}{2AC} \left[ \arctan\left(\frac{2u+B}{C}\right) + \arctan\left(\frac{2u-B}{C}\right) \right],
\end{equation}
where $A=\sqrt{1+10^{-10}}$, $B=\sqrt{2(A+1)}$, and $C=\sqrt{2(A-1)}$.

In all experiments, $F$ is the quartic potential $F(s)=\frac14(s^2-1)^2$ of Section~\ref{subsec:model} with the splitting \eqref{eq:semiconvex-splitting} for $\gamma=1$. The spatial domain is a rectangular domain in $\mathbb R^d$, $d\in\{2,3\}$, discretized by a uniform Cartesian mesh with $N_x\times N_y$ cells in 2D and $N_x\times N_y\times N_z$ cells in 3D, and $x_K$ is the center of $K$.
At each time step, we solve the nonlinear system by Newton's method and stop when the maximum norm of the Newton increment is below $10^{-8}$.
Theorem~\ref{thm:discrete-solvability} guarantees that a solution exists at each time step, and Proposition~\ref{prop:discrete-uniqueness} shows that it is unique for $M_0$.
In all experiments of this section, Newton's method met this stopping criterion at every time step.
We measure the cell-average error $e^n=\{e_K^n\}_{K\in\mesh}$, $e_K^n:=u_K^n-(P_{\mesh}u(\cdot,t_n))_K$, in the norms $\|e^n\|_{L^2(\Omega)}=(e^n,e^n)_{\mesh}^{1/2}$ and $\|e^n\|_{L^\infty(\Omega)}=\max_{K\in\mesh}|e_K^n|$ of the piecewise constant function $e^n$.

%-------------------------------------------------------------------------
\subsection[Accuracy and convergence]{{Accuracy and convergence}}
\label{subsec:exp_convergence}
%-------------------------------------------------------------------------

% \paragraph{problem setting}
We first examine the observed convergence orders of the scheme by the method of manufactured solutions on $\Omega=(0,1)^2$ with $T=1$ and, for all mobility models, $\eps=0.1$.
With a source term $f$ added to the right-hand side of the first equation in \eqref{eq:model}, the function $u(x,y,t)=\frac{1}{2}t^3\cos(\pi x)\cos(\pi y)$ is an exact solution, where
\[
f(x,y,t)
=\frac{3t^{3-\alpha}}{\Gamma(4-\alpha)}
\cos(\pi x)\cos(\pi y)
+M(u)
\bigl(2\eps^2\pi^2u+F'(u)\bigr).
\]
In the scheme, we replace \eqref{eq:scheme-evolution} by $\bar q_K^n\bigl((\dtaualpha u_K)^n-f(x_K,t_n)\bigr)=-\mu_K^n$.

% \paragraph{spatial convergence test}
For $\alpha=0.5$, we refine the Cartesian mesh with $N_x=N_y$ and take $N=N_x^2$ time steps at each resolution, so that $\tau=T/N=N_x^{-2}$.
Table~\ref{tab:spatial_error} shows that, for all three mobility models, the observed orders approach $2$, so the errors decrease like $\mathcal{O}(N_x^{-2})$ along this coupled refinement path.

\begin{table}[htbp]
    \centering
    \resizebox{\textwidth}{!}{
    \begin{tabular}{c|cccc|cccc|cccc}
        \hline
        \multirow{2}{*}{$N_x$} & \multicolumn{4}{c|}{$M_0$} & \multicolumn{4}{c|}{$M_1$} & \multicolumn{4}{c}{$M_2$} \\
        \cline{2-13}
        & $L^2$ err & Order & $L^\infty$ err & Order & $L^2$ err & Order & $L^\infty$ err & Order & $L^2$ err & Order & $L^\infty$ err & Order \\
        \hline
        $  4$ & 1.48e-2 & -    & 2.35e-2 & -    & 1.27e-2 & -    & 1.88e-2 & -    & 8.20e-3 & -    & 9.79e-3 & -    \\
        $  8$ & 5.46e-3 & 1.44 & 9.17e-3 & 1.35 & 4.68e-3 & 1.44 & 7.09e-3 & 1.41 & 2.96e-3 & 1.47 & 4.10e-3 & 1.26 \\
        $ 16$ & 1.54e-3 & 1.82 & 2.63e-3 & 1.80 & 1.32e-3 & 1.82 & 2.02e-3 & 1.81 & 8.42e-4 & 1.81 & 1.20e-3 & 1.78 \\
        $ 32$ & 4.04e-4 & 1.93 & 6.91e-4 & 1.93 & 3.47e-4 & 1.93 & 5.32e-4 & 1.93 & 2.22e-4 & 1.92 & 3.15e-4 & 1.93 \\
        $ 64$ & 1.03e-4 & 1.97 & 1.76e-4 & 1.97 & 8.86e-5 & 1.97 & 1.36e-4 & 1.97 & 5.70e-5 & 1.96 & 8.06e-5 & 1.97 \\
        $ 128$ & 2.60e-5 & 1.99 & 4.45e-5 & 1.99 & 2.24e-5 & 1.99 & 3.43e-5 & 1.98 & 1.44e-5 & 1.98 & 2.04e-5 & 1.98 \\
        \hline
    \end{tabular}
    }
    \caption{Discrete $L^2$ and $L^\infty$ errors at $T=1$ and observed orders for $\alpha=0.5$, with $N=N_x^2$ time steps.}
    \label{tab:spatial_error}
\end{table}

Next, we investigate the temporal convergence for the fractional orders $\alpha\in\{0.4,0.6,0.8\}$ on the fixed spatial mesh with $N_x = N_y = 128$, varying the time step $\tau = T/N$.
As summarized in Table~\ref{tab:temporal_error}, the observed orders approach $1$ for all three mobility models.
This suggests that the temporal error is dominated not by the L1 approximation, whose truncation error is $\mathcal{O}(\tau^{2-\alpha})$ for smooth solutions (see \cite[Section~3]{2007LinXu}), but by the first-order consistency error of the splitting \eqref{eq:semiconvex-splitting}, namely the term $\gamma(u^n-u^{n-1})$, with a smaller first-order contribution from the secant average for $M_1$ and $M_2$.
Both ingredients enter the proof of the energy bound in Theorem~\ref{thm:discrete-energy-stability}, which holds without a time-step restriction.

\begin{table}[htbp]
    \centering
    \resizebox{\textwidth}{!}{
    \begin{tabular}{c c | c c c c | c c c c | c c c c}
        \hline
        \multirow{2}{*}{$\alpha$} & \multirow{2}{*}{$N$} & \multicolumn{4}{c|}{$M_0$} & \multicolumn{4}{c|}{$M_1$} & \multicolumn{4}{c}{$M_2$} \\
        \cline{3-14}
        & & $L^2$ err & Order & $L^\infty$ err & Order & $L^2$ err & Order & $L^\infty$ err & Order & $L^2$ err & Order & $L^\infty$ err & Order \\
        \hline
        \multirow{5}{*}{0.4} 
        & $  50$ & 1.40e-2 & -    & 2.50e-2 & -    & 1.26e-2 & -    & 2.12e-2 & -    & 9.38e-3 & -    & 1.34e-2 & -    \\
        & $ 100$ & 7.34e-3 & 0.93 & 1.30e-2 & 0.94 & 6.57e-3 & 0.94 & 1.10e-2 & 0.95 & 4.83e-3 & 0.96 & 6.88e-3 & 0.97 \\
        & $ 200$ & 3.77e-3 & 0.96 & 6.67e-3 & 0.97 & 3.36e-3 & 0.97 & 5.58e-3 & 0.98 & 2.45e-3 & 0.98 & 3.48e-3 & 0.98 \\
        & $ 400$ & 1.91e-3 & 0.98 & 3.37e-3 & 0.99 & 1.69e-3 & 0.99 & 2.81e-3 & 0.99 & 1.23e-3 & 0.99 & 1.75e-3 & 0.99 \\
        & $ 800$ & 9.54e-4 & 1.00 & 1.68e-3 & 1.00 & 8.47e-4 & 1.00 & 1.40e-3 & 1.00 & 6.15e-4 & 1.01 & 8.69e-4 & 1.01 \\
        \hline
        \multirow{5}{*}{0.6} 
        & $  50$ & 8.80e-3 & -    & 1.63e-2 & -    & 8.07e-3 & -    & 1.42e-2 & -    & 6.32e-3 & -    & 9.31e-3 & -    \\
        & $ 100$ & 4.74e-3 & 0.89 & 8.77e-3 & 0.90 & 4.34e-3 & 0.89 & 7.60e-3 & 0.90 & 3.40e-3 & 0.90 & 4.99e-3 & 0.90 \\
        & $ 200$ & 2.48e-3 & 0.93 & 4.59e-3 & 0.94 & 2.27e-3 & 0.93 & 3.97e-3 & 0.94 & 1.78e-3 & 0.93 & 2.61e-3 & 0.93 \\
        & $ 400$ & 1.28e-3 & 0.96 & 2.35e-3 & 0.96 & 1.17e-3 & 0.96 & 2.04e-3 & 0.96 & 9.14e-4 & 0.96 & 1.34e-3 & 0.96 \\
        & $ 800$ & 6.46e-4 & 0.98 & 1.19e-3 & 0.98 & 5.91e-4 & 0.98 & 1.03e-3 & 0.99 & 4.62e-4 & 0.98 & 6.78e-4 & 0.99 \\
        \hline
        \multirow{5}{*}{0.8} 
        & $  50$ & 3.61e-3 & -    & 6.88e-3 & -    & 3.24e-3 & -    & 5.73e-3 & -    & 2.35e-3 & -    & 3.34e-3 & -    \\
        & $ 100$ & 2.11e-3 & 0.77 & 4.03e-3 & 0.77 & 1.92e-3 & 0.76 & 3.41e-3 & 0.75 & 1.44e-3 & 0.71 & 2.06e-3 & 0.69 \\
        & $ 200$ & 1.18e-3 & 0.84 & 2.25e-3 & 0.84 & 1.08e-3 & 0.83 & 1.93e-3 & 0.82 & 8.27e-4 & 0.80 & 1.20e-3 & 0.78 \\
        & $ 400$ & 6.40e-4 & 0.88 & 1.22e-3 & 0.89 & 5.87e-4 & 0.88 & 1.05e-3 & 0.87 & 4.56e-4 & 0.86 & 6.69e-4 & 0.84 \\
        & $ 800$ & 3.37e-4 & 0.92 & 6.41e-4 & 0.93 & 3.10e-4 & 0.92 & 5.55e-4 & 0.92 & 2.43e-4 & 0.91 & 3.58e-4 & 0.90 \\
        \hline
    \end{tabular}
    }\caption{Discrete $L^2$ and $L^\infty$ errors at $T=1$ and observed orders for $\alpha\in\{0.4,0.6,0.8\}$, with $N_x=N_y=128$.}
    \label{tab:temporal_error}
\end{table}

%-------------------------------------------------------------------------
\subsection[Spinodal decomposition]{{Spinodal decomposition}}
\label{subsec:exp_spinodal}
%-------------------------------------------------------------------------

% \paragraph{problem setting}
{\emergencystretch=2em
To examine the structure-preserving properties in a strongly nonlinear regime, we simulate spinodal decomposition on $\Omega=(0,2\pi)^2$ and $\Omega=(0,2\pi)^3$.
In both cases, the final time is $T=100$, the interfacial width parameter is $\eps=0.05$, and the fractional order is $\alpha=0.5$.
The initial data are a small random perturbation of the unstable mixed state $u=0$: we set $u_K^0=0.01\,\xi_K$ for every $K\in\mesh$, where the numbers $\xi_K$ are independent and uniformly distributed in $[-1,1]$.
\par}

Figure~\ref{fig:spinodal_snapshots} shows the 2D solution at $t=10$, $15$, $60$, and $100$.
Starting from the random initial data, the phases separate rapidly and the domains merge over time.
The three mobilities produce nearly the same domain pattern. For the regularized degenerate mobility $M_2$, the solution approaches the values $\pm1$ much more slowly than for $M_0$ and $M_1$ (see also Figure~\ref{fig:str-prev2}(a)), which is consistent with the small values of $M_2$ near $u=\pm1$.

\begin{figure}[htbp]
    \centering
    \includegraphics[width=0.95\textwidth]{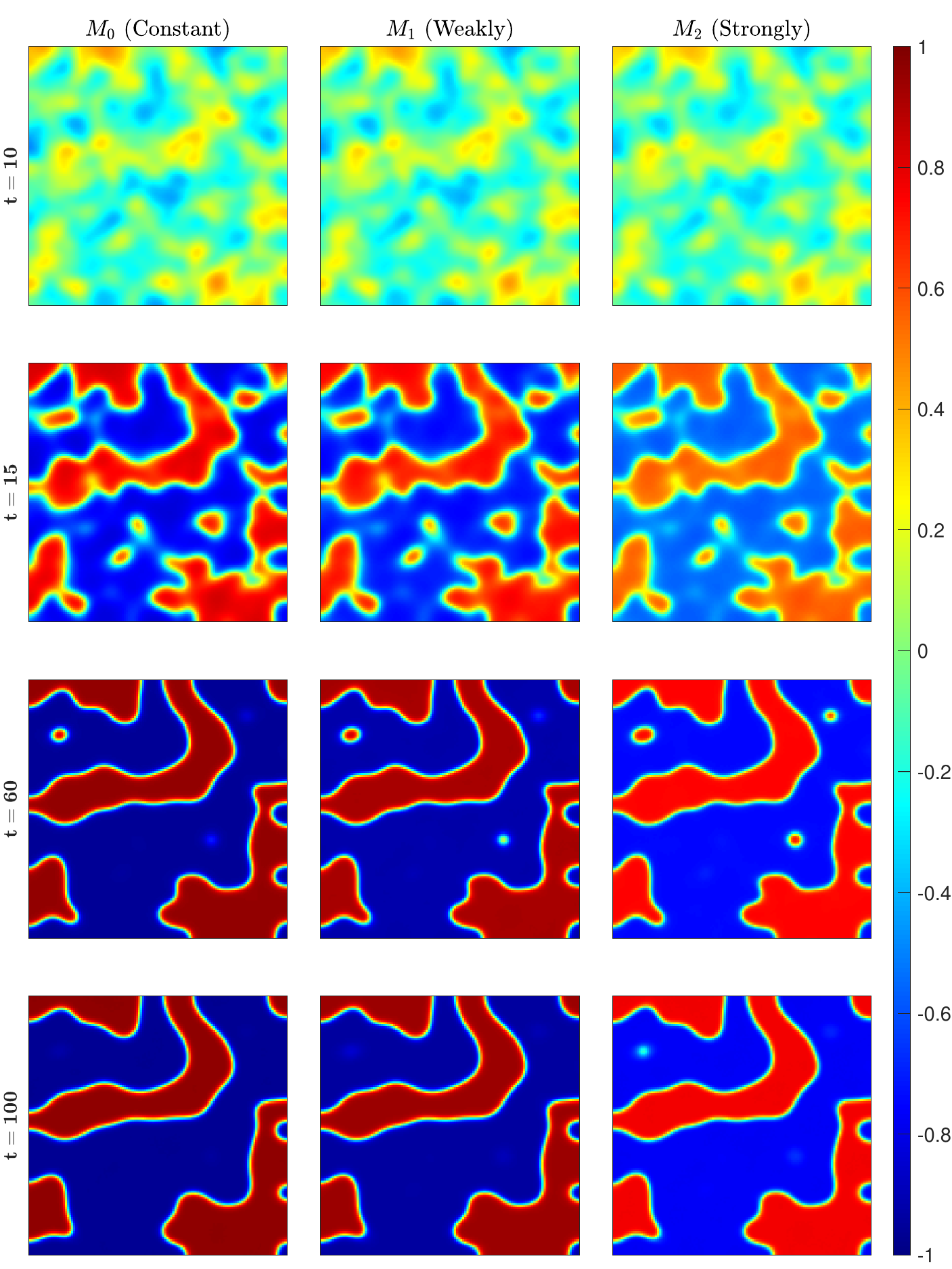}
    \caption{Snapshots of the 2D spinodal decomposition at $t=10$, $15$, $60$, and $100$ (rows) for the three mobilities: (left column) $M_0$; (center column) $M_1$; (right column) $M_2$. {(Simulation parameters: $N_x=N_y=128$, $N=200$)}}
    \label{fig:spinodal_snapshots}
\end{figure}

For the 2D and 3D simulations and all three mobility models, panels (a) of Figures~\ref{fig:str-prev2} and~\ref{fig:str-prev3} show that $\max_{K\in\mesh}|u_K^n|\le1$ at every time level, in agreement with the maximum bound principle in Theorem~\ref{thm:bound-preservation}.
Panels (b) show that the discrete energy $\energy_{\mesh}(u^n)$ stays below its initial value $\energy_{\mesh}(u^0)$, in agreement with Theorem~\ref{thm:discrete-energy-stability}.
Panels (c) show that the accumulated pairing of Section~\ref{sec:energy-stability} at $t_m$, $1\le m\le N$, stays below the initial discrete energy. The proof of Theorem~\ref{thm:discrete-energy-stability} bounds it by $\energy_{\mesh}(u^0)-\energy_{\mesh}(u^m)$, which is at most $\energy_{\mesh}(u^0)$ because the quartic potential is nonnegative.

\begin{figure}[htbp]
    \centering
    \includegraphics[width=0.98\textwidth]{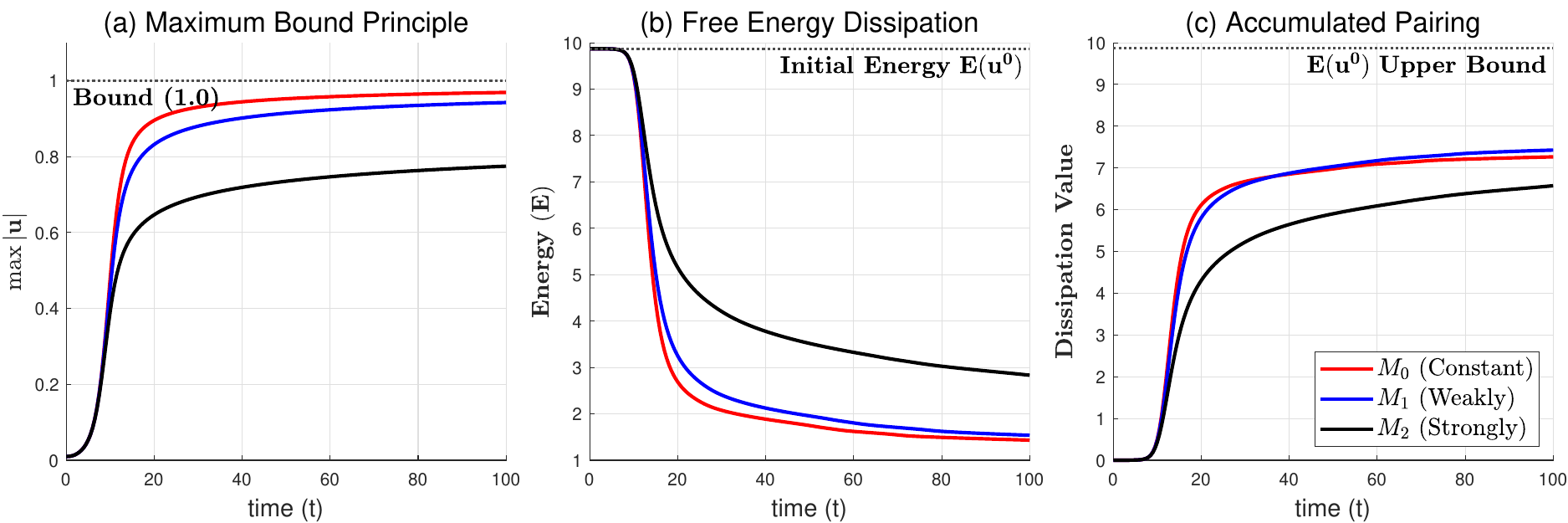}
    \caption{Structure-preserving properties in the 2D spinodal decomposition: (a) $\max_{K\in\mesh}|u_K^n|$; (b) the discrete energy $\energy_{\mesh}(u^n)$; (c) the accumulated pairing. Dotted lines mark the bound $1$ in (a) and the initial discrete energy in (b) and (c). {(Simulation parameters: $N_x=N_y=128$, $N=200$)}}
    \label{fig:str-prev2}
\end{figure}
\begin{figure}[htbp]
    \centering
    \includegraphics[width=0.98\textwidth]{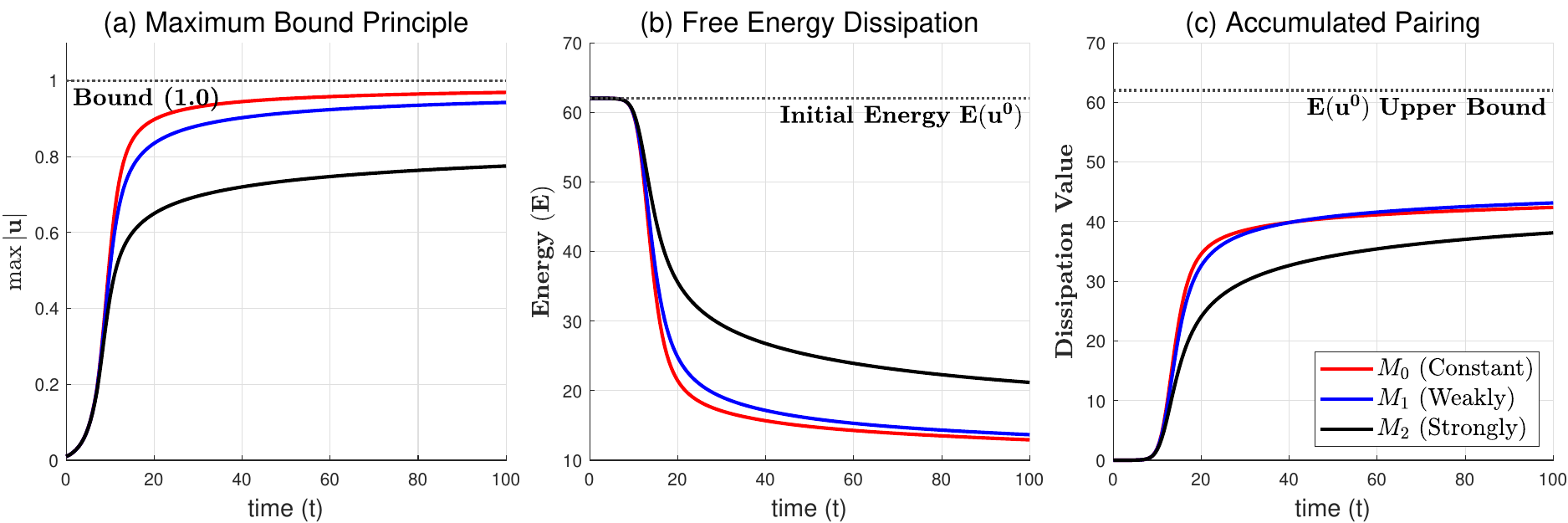}
    \caption{Structure-preserving properties in the 3D spinodal decomposition: (a) $\max_{K\in\mesh}|u_K^n|$; (b) the discrete energy $\energy_{\mesh}(u^n)$; (c) the accumulated pairing. Dotted lines mark the bound $1$ in (a) and the initial discrete energy in (b) and (c). {(Simulation parameters: $N_x=N_y=N_z=64$, $N=200$)}}
    \label{fig:str-prev3}
\end{figure}

%-------------------------------------------------------------------------
\subsection[Relaxation of a star-shaped droplet]{{Relaxation of a star-shaped droplet}}
\label{subsec:exp_mobility_comparison}
%-------------------------------------------------------------------------

In the final experiment, we compare the three mobility models in the relaxation of a star-shaped droplet in 3D and monitor the structure-preserving properties of the scheme.
The computation uses $\Omega=(0,1)^3$, the interfacial width parameter $\eps=0.02$, the fractional order $\alpha=0.5$, the final time $T=400$, and $N_x=N_y=N_z=64$.
The initial condition is
\[
u_0(x)=-\tanh\left(\frac{r(x)-R_{\text{star}}(\theta,\phi)}{\sqrt{2}\,\eps}\right),
\]
where $r(x)=|x-x_c|$ is the distance from the center $x_c=(0.5,0.5,0.5)$, $\theta$ and $\phi$ are the azimuthal and polar angles of $x-x_c$, and $R_{\text{star}}(\theta,\phi)=0.25+0.1\sin(4\theta)\sin(4\phi)$.

Figure~\ref{fig:mobility_compare}(a)--(c) show the zero-level isosurface $\{u=0\}$ for $M_2$ at $t=0$, $100$, and $400$.
To quantify the shape, we measure the minimum and maximum distances from $x_c$ to this isosurface (see Figure~\ref{fig:mobility_compare}(d)).
As shown in Figure~\ref{fig:mobility_compare}(e), for $M_0$ and $M_1$, the two distances converge as the surface becomes smoother, which indicates a transition to a nearly spherical shape. The droplet then shrinks and vanishes, and the distances drop to zero.
For the regularized degenerate mobility $M_2$, the relaxation is much slower: at $t=400$, the two distances have not yet merged.
\begin{figure}[htbp]
    \centering
    \includegraphics[width=0.98\textwidth]{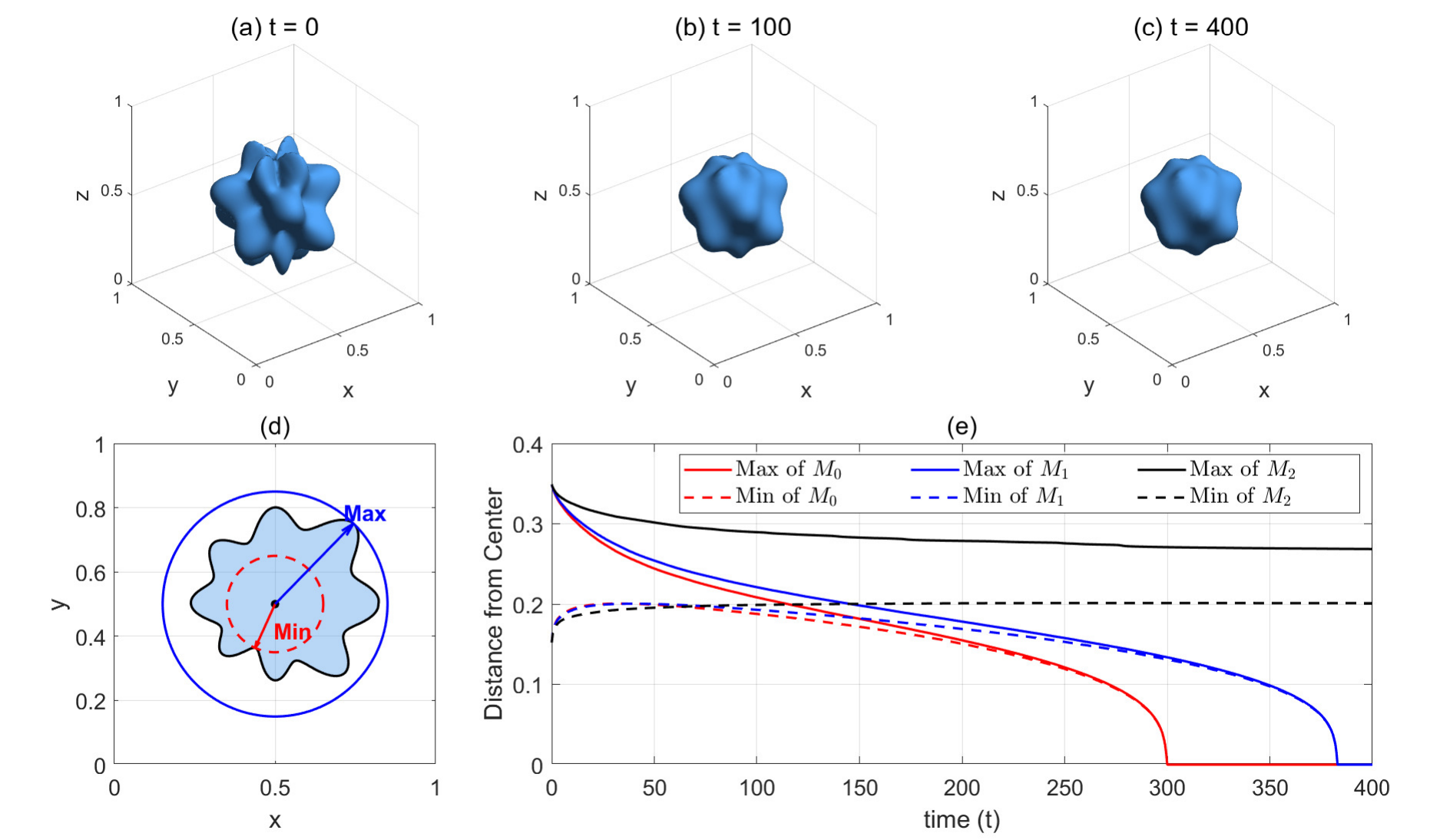}
    \caption{Relaxation of the 3D star-shaped droplet: (a)--(c) the zero-level isosurface for $M_2$ at $t=0$, $100$, and $400$; (d) a two-dimensional sketch of the minimum and maximum distances from the center; (e) these distances over time for the three mobility models. {(Simulation parameters: $N_x=N_y=N_z=64$, $N=800$)}}
    \label{fig:mobility_compare}
\end{figure}

Figure~\ref{fig:mobility_compare1} shows the quantities of Figures~\ref{fig:str-prev2} and~\ref{fig:str-prev3} for this experiment.
Away from the droplet, the cell values equal $-1$ to machine precision, so $\max_{K\in\mesh}|u_K^n|$ equals $1$ to machine precision at every time level. Panel (a) therefore shows only that no cell value exceeds $1$ in absolute value, in agreement with Theorem~\ref{thm:bound-preservation}.
For all three mobility models, the discrete energy stays below its initial value, and the accumulated pairing stays below the initial discrete energy, in agreement with Theorem~\ref{thm:discrete-energy-stability}.

\begin{figure}[htbp]
    \centering
    \includegraphics[width=0.98\textwidth]{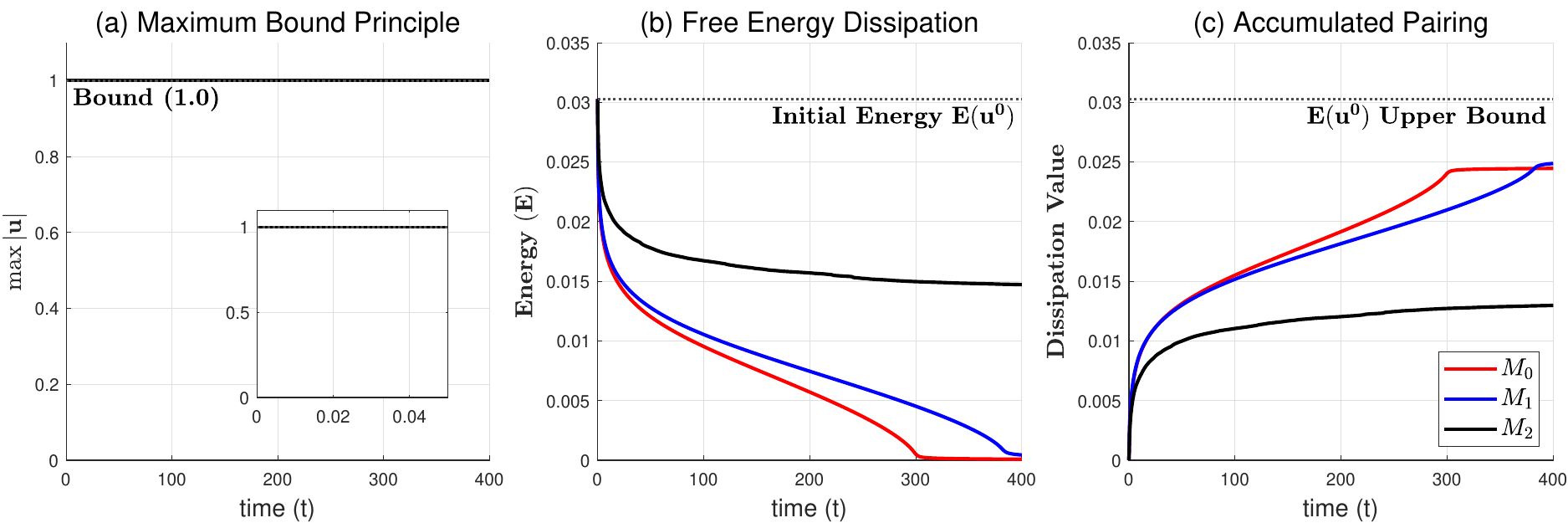}
    \caption{Structure-preserving properties in the 3D star-shaped droplet experiment: (a) $\max_{K\in\mesh}|u_K^n|$; (b) the discrete energy $\energy_{\mesh}(u^n)$; (c) the accumulated pairing. Dotted lines mark the bound $1$ in (a) and the initial discrete energy in (b) and (c). The inset in (a) magnifies $0\le t\le0.05$. {(Simulation parameters: $N_x=N_y=N_z=64$, $N=800$)}}
    \label{fig:mobility_compare1}
\end{figure}

}

\clearpage
\section[Conclusion and future work]{{Conclusion and future work}}
\label{sec:conc}
We have developed a structure-preserving finite volume scheme for the time-fractional Allen--Cahn equation with nondegenerate state-dependent mobility, combining the two-point flux approximation of diffusion, the uniform-grid L1 approximation of the Caputo derivative, and a semiconvex splitting of the potential with a secant average of the inverse mobility that multiplies the entire L1 derivative.
On every admissible finite volume mesh and uniform time grid, the scheme has at least one solution, and every solution satisfies the maximum bound principle and an energy bound relative to the initial energy. All three results hold without a time-step restriction.
Under an additional condition on the variation of the inverse mobility, the fully discrete solution is unique.
The energy bound and the time regularity required for compactness follow from the key estimate, a lower bound for the accumulated pairing in terms of the accumulated L1 quadratic form.
Under a mesh regularity condition with a constant independent of the mesh size, any sequence of fully discrete solutions has a subsequence that converges strongly in $L^2(\OmegaT)$ to a weak solution of \eqref{eq:model}--\eqref{eq:initial-condition} as the mesh size and the time step tend to zero.
Consequently, a weak solution of \eqref{eq:model}--\eqref{eq:initial-condition} exists in one and two space dimensions, and in three dimensions whenever a sequence of meshes satisfying this condition exists.

In the tests of Section~\ref{subsec:exp_convergence} with a manufactured solution that is smooth in time, where $N_x$ is the number of cells in each coordinate direction, the errors decrease like $\mathcal{O}(N_x^{-2})$ along the coupled refinement path $\tau=N_x^{-2}$, and the observed orders in time approach $1$ on the fixed mesh with $N_x=128$.
The two- and three-dimensional experiments in Sections~\ref{subsec:exp_spinodal} and~\ref{subsec:exp_mobility_comparison} illustrate the maximum bound principle, the energy bound relative to the initial energy, and the upper bound for the accumulated pairing obtained in the proof of Theorem~\ref{thm:discrete-energy-stability}.

Solutions of time-fractional equations may have singular time derivatives near $t=0$, even for regular initial data.
Quantifying the influence of this singularity on numerical accuracy is a natural extension of the present analysis, and constructing schemes that treat this singularity accurately and still satisfy the maximum bound principle and the energy bound is an important research direction, which we are pursuing in follow-up work.
By Remark~\ref{rem:nonuniform-time-grids}, both properties already hold on nonuniform time grids such as graded grids, and this direction includes extending the convergence analysis of Sections~\ref{subsec:uniform-estimates-compactness} and~\ref{sec:convergence-analysis} to such grids.
Another direction is the question left open in Remark~\ref{rem:secant-coefficient-choice}: an energy estimate without a time-step restriction for schemes in which the inverse mobility is evaluated at the previous time level.

\paragraph{Funding}
Y.J. was supported by the National Research Foundation of Korea (NRF) grant
funded by the Korea government (MSIT) (No.~RS-2026-25479175).
K.K. was supported by the National Research Foundation of Korea (NRF) grants
funded by the Korea government (MSIT) (Nos.~00560003 and 2023R1A2C1003120).

\clearpage
\appendix
\section[Proofs of auxiliary results]{{Proofs of auxiliary results}}
\label{app:technical-estimates}

\subsection[Second differences of the L1 weights]{{Second differences of the L1 weights}}
\label{app:l1-weights}

Throughout this appendix, $\beta:=1-\alpha$.

\begin{lemma}[{Second differences of the L1 weights}]
\label{lem:uniform-l1-weight-properties}
The L1 weights \eqref{eq:l1-weights} satisfy, for every integer $r\ge2$,
\begin{equation}
  a_{r-1}-2a_r+a_{r+1}
  \ge\alpha\beta(1+\alpha)2^{-2-\alpha}r^{-2-\alpha}.
  \label{eq:l1-second-difference-lower-bound}
\end{equation}
\end{lemma}

\begin{proof}
Repeated integration with $g(x)=\beta x^{-\alpha}$ gives, for $r\ge2$,
\begin{equation*}
  a_{r-1}-2a_r+a_{r+1}
  =\int_{[0,1]^3}g''(r-1+s+t+u)\,\dd s\,\dd t\,\dd u,
\end{equation*}
where $g''(x)=\alpha\beta(1+\alpha)x^{-2-\alpha}$. Since $g''$ is positive
and decreasing and $0<r-1+s+t+u\le2r$, this gives
\eqref{eq:l1-second-difference-lower-bound}.
\end{proof}

\subsection[L1 matrix entries and edge weights]{{L1 matrix entries and edge weights}}
\label{app:l1-matrix}

{The edge weights \eqref{eq:l1-symmetric-edge-weight} are
$\omega_{01}=a_0=1$ for $m=1$ and, for $m\ge2$, $0\le i<j\le m$, and
$r:=j-i$,}
\begin{equation}
  \omega_{ij}
  =
  \begin{cases}
    \dfrac12(2a_0-a_1),
    & r=1,\quad i=0\ \text{or}\ j=m,\\[2mm]
    \dfrac12(2a_0-2a_1+a_2),
    & r=1,\quad 1\le i<j\le m-1,\\[2mm]
    \dfrac12(a_{r-1}-a_r),
    & r\ge2,\quad
      (i=0,\ j<m)\ \text{or}\ (i>0,\ j=m),\\[2mm]
    \dfrac12(a_{r-1}-2a_r+a_{r+1}),
    & r\ge2,\quad 1\le i<j\le m-1,\\[2mm]
    \dfrac12a_{m-1},
    & (i,j)=(0,m).
  \end{cases}
  \label{eq:l1-symmetric-edge-weight-cases}
\end{equation}

\begin{proof}[Proof of Lemma~\ref{lem:uniform-l1-matrix-structure}]
{
{
For part (i), to write the
entries of $A_m$ uniformly for $0\le i,j\le m$, we extend $L_m$ by zero
outside its original index range. Using $A_m=D^\top L_mD$ and
$(Dx)_p=x_p-x_{p-1}$, we obtain
\begin{equation*}
  (A_m)_{ij}
  =
  (L_m)_{ij}
  -(L_m)_{i,j+1}
  -(L_m)_{i+1,j}
  +(L_m)_{i+1,j+1}.
\end{equation*}
{%
Substituting the entries of $L_m$, we obtain, for $i\ne j$,
\begin{equation}
  (A_m)_{ij}=
  \begin{cases}
  -a_0,&j=i+1,\quad0\le i\le m-1,\\
  0,&j\ge i+2,\\
  a_i-a_{i-1},&j=0,\quad1\le i\le m-1,\\
  -\bigl(a_{i-j-1}-2a_{i-j}+a_{i-j+1}\bigr),&1\le j<i\le m-1,\\
  -a_{m-1},&(i,j)=(m,0),\\
  a_{m-j}-a_{m-j-1},&i=m,\quad1\le j\le m-1.
  \end{cases}
  \label{eq:l1-off-diagonal-entry-cases}
\end{equation}
The positivity of the L1 weights gives $-a_0,-a_{m-1}\le0$.
Their monotonicity gives $a_i-a_{i-1}\le0$ and
$a_{m-j}-a_{m-j-1}\le0$, and their discrete convexity gives
$a_{i-j-1}-2a_{i-j}+a_{i-j+1}\ge0$ (see
\eqref{eq:l1-weight-properties}). Hence $(A_m)_{ij}\le0$ for every $i\ne j$.
}
Since $D\mathbf1=0$, both $A_m\mathbf1$ and $A_m^\top\mathbf1$ vanish.
}

{
For part (ii), set $\mathbf v=(v^0,\ldots,v^m)^\top$, so that
$D\mathbf v=(\delta v^1,\ldots,\delta v^m)^\top$.
Part (i) gives $\omega_{ij}\ge0$ and $S_m\mathbf1=0$.
{Using $A_m=D^\top L_mD$, symmetry of the inner product,
and the zero row sums of $S_m$, we obtain}
\begin{equation*}
  \begin{aligned}
  \sum_{n=1}^m\langle\delta v^n,(L_m\delta v)_n\rangle
  &=\sum_{i,j=0}^m(S_m)_{ij}\langle v^i,v^j\rangle\\
  &=\sum_{0\le i<j\le m}\omega_{ij}\|v^j-v^i\|^2.
  \end{aligned}
\end{equation*}
This proves \eqref{eq:l1-quadratic-edge-representation}.
}

{
Finally, we verify \eqref{eq:l1-symmetric-edge-weight-cases}. If $m=1$, a
direct calculation gives
\begin{equation*}
  A_1=S_1=
  \begin{pmatrix}
    1&-1\\
    -1&1
  \end{pmatrix},
\end{equation*}
and hence $\omega_{01}=1$. For $m\ge2$, substituting
\eqref{eq:l1-off-diagonal-entry-cases} for $(A_m)_{ij}$ and $(A_m)_{ji}$ into
\eqref{eq:l1-symmetric-edge-weight} yields
\eqref{eq:l1-symmetric-edge-weight-cases}.
}
}
\end{proof}

\subsection[Projection onto cell averages]{{Projection onto cell averages}}
\label{app:initial-projection}

\begin{proof}[Proof of Lemma~\ref{lem:initial-projection}]
{
{
Let $v_K:=(P_\ell v)_K$, let $v_\sigma$ be the average of the trace of $v$ on
$\sigma$, and set $h_K:=\operatorname{diam}(K)$ and $w_K:=v-v_K$. For an
interior face $\sigma=K|L$, the points $x_K\in K$ and $x_L\in L$ lie on
opposite sides of the hyperplane containing $\sigma$, so
$d_\sigma=d_{K,\sigma}+d_{L,\sigma}$ by the orthogonality condition of
Section~\ref{subsec:finite-volume-discretization}. The weighted Cauchy inequality gives
\begin{equation*}
  \frac{(v_L-v_K)^2}{d_\sigma}
  \le\frac{(v_\sigma-v_K)^2}{d_{K,\sigma}}
  +\frac{(v_L-v_\sigma)^2}{d_{L,\sigma}}.
\end{equation*}
Summing over interior faces, applying Jensen's inequality and
$d_{K,\sigma}\ge\vartheta h_K$ there, and adding the nonnegative contributions
of the boundary faces to the last sum, we obtain
\begin{equation*}
  \begin{aligned}
  |P_\ell v|_{1,\mesh_\ell}^2
  &\le\sum_{K\in\mesh_\ell}
  \sum_{\sigma\in\mathcal E_{K,\mathrm{int}}}
  \frac{|\sigma|}{d_{K,\sigma}}(v_\sigma-v_K)^2\\
  &\le\vartheta^{-2}\sum_{K\in\mesh_\ell}h_K^{-2}
  \sum_{\sigma\in\mathcal E_K}
  d_{K,\sigma}\int_\sigma|w_K|^2\,\dd S.
  \end{aligned}
\end{equation*}
Since $K$ is convex and $x_K\in K$, we have
$(x-x_K)\cdot n_{K,\sigma}=d_{K,\sigma}$ for every $\sigma\in\mathcal E_K$
and $x\in\sigma$. Applying the Gauss--Green
formula to $(x-x_K)w_K^2$, we therefore obtain
\begin{equation*}
  \sum_{\sigma\in\mathcal E_K}d_{K,\sigma}\int_\sigma|w_K|^2\,\dd S
  =d\|w_K\|_{L^2(K)}^2
  +2\int_Kw_K(x-x_K)\cdot\nabla v\,\dd x.
\end{equation*}
Since $|x-x_K|\le h_K$, the Poincar\'e inequality
$\|w_K\|_{L^2(K)}\le(h_K/\pi)\|\nabla v\|_{L^2(K)}$
\cite[equations~(1.4) and~(1.9)]{1960Payne}
and \cite[Theorem~3.2]{2003Bebendorf} gives
\begin{equation*}
  \sum_{\sigma\in\mathcal E_K}d_{K,\sigma}\int_\sigma|w_K|^2\,\dd S
  \le h_K^2\left(\frac d{\pi^2}+\frac2\pi\right)
  \|\nabla v\|_{L^2(K)}^2.
\end{equation*}
Substitution into the bound for $|P_\ell v|_{1,\mesh_\ell}^2$
proves the first assertion of
Lemma~\ref{lem:initial-projection}. See
\cite[Lemma~9.4]{2000EymardGallouetHerbin} for a related estimate in two
dimensions.
}

Cell averaging preserves the bound $|u_0|\le1$, and
{applying the same Poincar\'e inequality on each cell
and summing over the mesh gives}
\begin{equation*}
  \|P_\ell u_0-u_0\|_{L^2(\Omega)}^2
  \le\frac{h_\ell^2}{\pi^2}\|\nabla u_0\|_{L^2(\Omega)}^2\longrightarrow0.
\end{equation*}
This proves \eqref{eq:initial-projection-properties}. Finally, we apply
the first assertion of Lemma~\ref{lem:initial-projection}
with $v=u_0$ and use
$F(u_{\ell,K}^0)\le F_{\max}$ in the definition of the discrete energy to
obtain the stated bound for the initial energies.
}
\end{proof}

\subsection[Lower bounds for the L1 edge weights]{{Lower bounds for the L1 edge weights}}
\label{app:l1-coercivity}

\begin{proof}[Proof of Lemma~\ref{lem:l1-fractional-coercivity}]
{
Since $0<\beta<1$ and $c_\alpha/\beta=\alpha(1+\alpha)2^{-3-\alpha}<1/4$, we
have $c_\alpha<1/2$.
}%
{%
If $m=1$, Appendix~\ref{app:l1-matrix} gives
$\omega_{01}=1\ge c_\alpha$. Let $m\ge2$, fix $0\le i<j\le m$, and set
$r=j-i$. For $r=1$, the two possible coefficients in
\eqref{eq:l1-symmetric-edge-weight-cases} satisfy
\begin{equation*}
  \begin{aligned}
  \frac12(2a_0-a_1)
  &=\frac12\bigl(1+a_0-a_1\bigr)\ge\frac12,\\
  \frac12(2a_0-2a_1+a_2)
  &=\frac12\bigl(1+a_0-2a_1+a_2\bigr)\ge\frac12,
  \end{aligned}
\end{equation*}
by $a_0=1$ and \eqref{eq:l1-weight-properties}. Hence
$\omega_{ij}\ge c_\alpha$.

For $r\ge2$, every coefficient in
\eqref{eq:l1-symmetric-edge-weight-cases} is at least the coefficient
$\frac12(a_{r-1}-2a_r+a_{r+1})$ of an edge with $1\le i<j\le m-1$. The
coefficient for exactly one of $i,j$ in $\{0,m\}$ exceeds it by
$\frac12(a_r-a_{r+1})$, and the coefficient for $(i,j)=(0,m)$, where $r=m$,
exceeds it by $\frac12(2a_m-a_{m+1})$. Both differences are positive by
\eqref{eq:l1-weight-properties}. Hence
\eqref{eq:l1-second-difference-lower-bound} gives
$\omega_{ij}\ge\frac12\alpha\beta(1+\alpha)2^{-2-\alpha}r^{-2-\alpha}
=c_\alpha r^{-2-\alpha}$.
}
{
{Thus, for every integer $m\ge1$,}
\begin{equation}
  \omega_{ij}
  \ge
  c_\alpha(j-i)^{-2-\alpha},
  \qquad
  0\le i<j\le m.
  \label{eq:l1-all-edge-coercivity}
\end{equation}
Substituting this bound into the edge representation
\eqref{eq:l1-quadratic-edge-representation}, we obtain
\eqref{eq:l1-fractional-coercivity}.\qedhere

}
\end{proof}

\subsection[The \texorpdfstring{$H^s$}{Hs} seminorm of piecewise affine interpolants]{{The $H^s$ seminorm of piecewise affine interpolants}}
\label{app:affine-interpolation}

\begin{proof}[Proof of Lemma~\ref{lem:affine-fractional-time-estimate}]
{
{Set $p=1+2s=2+\alpha$ and $I^n=(t_{n-1},t_n)$ for $1\le n\le N$. We split the double integral over
$I^i\times I^j$ into three sums: $S_{\mathrm{same}}$ for $i=j$,
$S_{\mathrm{adj}}$ for $|i-j|=1$, and $S_{\mathrm{far}}$ for $|i-j|\ge2$.
Each sum includes all ordered pairs satisfying its index condition. For
$N=1$, $S_{\mathrm{adj}}$ and $S_{\mathrm{far}}$ are empty sums and vanish.}

{
For $t,r\in I^n$, we have
$\widehat v(t)-\widehat v(r)=(t-r)\delta v^n/\tau$. Integrating
$|t-r|^{-\alpha}$ over $I^n\times I^n$, we obtain
\begin{equation*}
  S_{\mathrm{same}}
  =\frac{2\tau^{-\alpha}}{(1-\alpha)(2-\alpha)}
  \sum_{n=1}^N\|\delta v^n\|_{L^2(\Omega)}^2
  \le\frac{2}{(1-\alpha)(2-\alpha)}Q_\tau(v).
\end{equation*}
For $r\in I^n$ and $t\in I^{n+1}$, set $r=t_n-\tau x$ and
$t=t_n+\tau y$, with $x,y\in(0,1)$. Then
$\widehat v(t)-\widehat v(r)=x\,\delta v^n+y\,\delta v^{n+1}$ and $t-r=\tau(x+y)$.
Using
\begin{equation*}
  \|x\,\delta v^n+y\,\delta v^{n+1}\|_{L^2(\Omega)}^2
  \le2(x+y)^2\bigl(\|\delta v^n\|_{L^2(\Omega)}^2
  +\|\delta v^{n+1}\|_{L^2(\Omega)}^2\bigr),
\end{equation*}
and counting both integration orders, we obtain
\begin{equation*}
  \begin{aligned}
  S_{\mathrm{adj}}
  &\le4K_\alpha\tau^{-\alpha}\sum_{n=1}^{N-1}
  \bigl(\|\delta v^n\|_{L^2(\Omega)}^2+\|\delta v^{n+1}\|_{L^2(\Omega)}^2\bigr)\\
  &\le8K_\alpha Q_\tau(v),
  \end{aligned}
\end{equation*}
where
\begin{equation*}
  K_\alpha=\int_0^1\int_0^1(x+y)^{-\alpha}\,\dd x\,\dd y
  =\frac{2^{2-\alpha}-2}{(1-\alpha)(2-\alpha)}.
\end{equation*}
}%
{%
For $r\in I^m$ and $t\in I^n$ with $k=n-m\ge2$, write
$r=t_{m-1}+\tau y$ and $t=t_{n-1}+\tau x$, with $x,y\in(0,1)$. Then
\begin{equation*}
  \begin{aligned}
  \widehat v(t)-\widehat v(r)
  ={}&(1-x)(1-y)(v^{n-1}-v^{m-1})\\
  &+(1-x)y(v^{n-1}-v^m)\\
  &+x(1-y)(v^n-v^{m-1})\\
  &+xy(v^n-v^m).
  \end{aligned}
\end{equation*}
The four coefficients are nonnegative, sum to one, and each has integral
$1/4$ over $(0,1)^2$. By convexity of the squared norm and
$t-r\ge\tau(k-1)$,
\begin{equation*}
  \int_{I^n}\int_{I^m}
  \frac{\|\widehat v(t)-\widehat v(r)\|_{L^2(\Omega)}^2}{|t-r|^p}\,\dd r\,\dd t
  \le\frac{\tau^{-\alpha}}{4(k-1)^p}
  \sum_{a,b\in\{0,1\}}\|v^{n-1+a}-v^{m-1+b}\|_{L^2(\Omega)}^2.
\end{equation*}
For each fixed $a,b$, the map $(m,n)\mapsto(m-1+b,n-1+a)$ is injective, and
the index difference $k+a-b$ of the image pair satisfies
$1\le k+a-b\le k+1\le3(k-1)$.
Thus, summing over $m,n$, the four choices of $a,b$, and both integration
orders gives
\begin{equation*}
  S_{\mathrm{far}}
  \le2\cdot3^p\tau^{-\alpha}
  \sum_{0\le i<j\le N}\frac{\|v^j-v^i\|_{L^2(\Omega)}^2}{(j-i)^p}
  =2\cdot3^pQ_\tau(v).
\end{equation*}
}

{Adding the bounds for $S_{\mathrm{same}}$,
$S_{\mathrm{adj}}$, and $S_{\mathrm{far}}$, we obtain}
\begin{equation*}
  [\widehat v]_{H^s(0,T;L^2(\Omega))}^2
  =S_{\mathrm{same}}+S_{\mathrm{adj}}+S_{\mathrm{far}}
  \le C_{\mathrm{aff},\alpha}Q_\tau(v),
\end{equation*}
where
$C_{\mathrm{aff},\alpha}:=\frac{2}{(1-\alpha)(2-\alpha)}+8K_\alpha+2\cdot3^{2+\alpha}$.
}
\end{proof}

\subsection[Right-sided Caputo tests and weak limits of the L1 derivative]{{Right-sided Caputo tests and weak limits of the L1 derivative}}
\label{app:caputo-duality}

\begin{lemma}[Right-sided Caputo tests]
\label{lem:right-caputo-test-consistency}
Let $\eta\in C^2([0,T])$ vanish on $[T-\delta,T]$ for some $0<\delta\le T$,
and define $\widetilde\partial_{T-}^\alpha\eta$ on $[0,T]$ by the formula
preceding Definition~\ref{def:weak-caputo-derivative}.
On the uniform time grid of Section~\ref{subsec:l1-approximation}, let
\begin{equation*}
  R_\tau^j\eta
  :=\kappaalpha
    \sum_{n=j}^{N-1}a_{n-j}\bigl(\eta(t_n)-\eta(t_{n+1})\bigr),
  \qquad 1\le j\le N,
\end{equation*}
with the empty sum interpreted as zero. Let $R_\tau\eta(t)=R_\tau^j\eta$
for $t\in(t_{j-1},t_j]$. Then
\begin{equation}
  \|R_\tau\eta-\widetilde\partial_{T-}^\alpha\eta\|_{L^\infty(0,T)}
  \le2C_{\alpha,T}\tau\|\eta''\|_{L^\infty(0,T)},
  \label{eq:right-l1-test-first-order-consistency}
\end{equation}
where $C_{\alpha,T}:=T^{1-\alpha}/\Gamma(2-\alpha)$.
\end{lemma}

\begin{proof}
Extend $\eta$ by zero for $t\ge T$. Since $\eta$ vanishes near $T$, the
extension is $C^2$ and its derivatives have the same bounds. A change of
variables gives
\begin{equation*}
  (\widetilde\partial_{T-}^\alpha\eta)(t)=-\frac1{\Gamma(1-\alpha)}
  \int_0^T r^{-\alpha}\eta'(t+r)\,\dd r.
\end{equation*}
Let $\widehat\eta$ be the continuous piecewise affine interpolant of
$\eta(t_0),\ldots,\eta(t_N)$, extended by zero for $t\ge T$.
By the definition \eqref{eq:l1-weights} of $a_k$, we have
$\int_{t_n}^{t_{n+1}}(r-t_j)^{-\alpha}\,\dd r=\tau^{1-\alpha}a_{n-j}/(1-\alpha)$.
Since $\Gamma(2-\alpha)=(1-\alpha)\Gamma(1-\alpha)$, we obtain
\begin{equation*}
  R_\tau^j\eta
  =-\frac1{\Gamma(1-\alpha)}
  \int_0^T r^{-\alpha}\widehat\eta'(t_j+r)\,\dd r.
\end{equation*}
Let $t\in(t_{j-1},t_j]$ and $0<r<T$. By the mean value theorem,
$|\widehat\eta'(t_j+r)-\eta'(t_j+r)|\le\tau\|\eta''\|_{L^\infty(0,T)}$
whenever $t_j+r$ is not a grid point, and $|t_j-t|\le\tau$ gives
$|\eta'(t_j+r)-\eta'(t+r)|\le\tau\|\eta''\|_{L^\infty(0,T)}$.
Integrating $r^{-\alpha}$ over $(0,T)$ proves
\eqref{eq:right-l1-test-first-order-consistency}.
\end{proof}

{For weak Caputo derivatives in continuous time,} an identification of derivative
limits under strong convergence of the functions is given in
\cite[Proposition~3.5]{2018LiLiu}. The following lemma establishes an
identification for the L1 approximation under weak convergence of the
functions and initial values.

\begin{lemma}[Weak limits of the L1 derivative]
\label{lem:hilbert-l1-limit-identification}
For each $\ell$, let $N_\ell$,
$\tau_\ell=T/N_\ell$, and $t_\ell^n$ be as in
Assumption~\ref{ass:convergence-discretizations}, and let
$v_\ell^0,\ldots,v_\ell^{N_\ell}\in L^2(\Omega)$. Define $v_\ell^+$ as in
Definition~\ref{def:time-reconstructions}, and let
$\mathcal D_\ell^\alpha v_\ell(t):=(\partial_{\tau_\ell}^\alpha v_\ell)^n$ for
$t\in I_\ell^n$ and $1\le n\le N_\ell$, where the L1 approximation of this
$L^2(\Omega)$-valued sequence with time step $\tau_\ell$ is defined in
Section~\ref{sec:energy-stability}. For cellwise finite volume functions,
this is \eqref{eq:discrete-caputo-reconstruction}.
Suppose that
\begin{equation*}
  v_\ell^+\rightharpoonup v\text{ in }L^2(0,T;L^2(\Omega)),
  \qquad
  v_\ell^0\rightharpoonup v_0\text{ in }L^2(\Omega),
  \qquad
  \sup_\ell\|\mathcal D_\ell^\alpha v_\ell\|_{L^2(0,T;L^2(\Omega))}<\infty.
\end{equation*}
Then there is a unique $w\in L^2(0,T;L^2(\Omega))$ such that, for every
$\eta\in C_c^\infty((-\infty,T);\mathbb R)$,
\begin{equation}
  \int_0^T w(t)\eta(t)\,\dd t
  =\int_0^T(v(t)-v_0)
       (\widetilde\partial_{T-}^\alpha\eta)(t)\,\dd t
  \qquad\text{in }L^2(\Omega),
  \label{eq:identified-weak-caputo-pairing}
\end{equation}
where both integrals are Bochner integrals. Moreover,
$\mathcal D_\ell^\alpha v_\ell\rightharpoonup w$ in $L^2(0,T;L^2(\Omega))$.
Thus $w$ represents the weak Caputo derivative of $v$ associated with
$v_0$, defined by \eqref{eq:weak-caputo-derivative}.
\end{lemma}

\begin{proof}
Boundedness in $L^2(0,T;L^2(\Omega))$ gives a subsequence on which
$\mathcal D_\ell^\alpha v_\ell\rightharpoonup w$.
Let $\eta$ be a scalar test function as in the statement. Since its support
is a compact subset of $(-\infty,T)$, its restriction to $[0,T]$ satisfies
the assumptions of Lemma~\ref{lem:right-caputo-test-consistency}. For fixed
$\ell$, we omit this index and write $\tau$, $N$, $I^n$, $v^n$,
$\kappa_\alpha$, $v_\tau^+$, and $\mathcal D_\tau^\alpha v$ for $\tau_\ell$,
$N_\ell$, $I_\ell^n$, $v_\ell^n$, $\kappa_{\alpha,\ell}$, $v_\ell^+$, and
$\mathcal D_\ell^\alpha v_\ell$. Let $\eta^n=\eta(t_n)$ with $t_n=n\tau$,
so that $\eta^N=0$, and let $\eta_\tau$ be the function equal to $\eta^n$ on
$I^n$ for $1\le n\le N$.
With $R_\tau\eta$ from
Lemma~\ref{lem:right-caputo-test-consistency}, set
\begin{equation*}
  \lambda^j:=\tau\kappa_\alpha\sum_{n=j}^N a_{n-j}\eta^n
  \quad(1\le j\le N),
  \qquad \lambda^{N+1}:=0.
\end{equation*}
{Since $\eta^N=0$, direct subtraction gives
$\lambda^j-\lambda^{j+1}=\tau R_\tau^j\eta$ for $1\le j\le N$.}
Reordering the finite sums and applying Abel summation gives the exact
Bochner integral identity
\begin{equation}
  \begin{aligned}
  \int_0^T\mathcal D_\tau^\alpha v\,\eta_\tau\,\dd t
  &=\sum_{j=1}^N(v^j-v^{j-1})\lambda^j\\
  &=\sum_{j=1}^N(v^j-v^0)(\lambda^j-\lambda^{j+1})\\
  &=\int_0^T(v_\tau^+-v^0)R_\tau\eta\,\dd t
  \qquad\text{in }L^2(\Omega).
  \end{aligned}
  \label{eq:discrete-caputo-integral-pairing}
\end{equation}

{
The smoothness of $\eta$ gives $\eta_{\tau_\ell}\to\eta$
uniformly on $(0,T]$, and the estimate
\eqref{eq:right-l1-test-first-order-consistency} gives
$R_{\tau_\ell}\eta\to\widetilde\partial_{T-}^{\alpha}\eta$ uniformly
on $(0,T)$.
Since $T<\infty$, for every fixed $z\in L^2(\Omega)$,
\begin{equation*}
  \eta_{\tau_\ell}z\to\eta z,
  \qquad
  (R_{\tau_\ell}\eta)z\to
  (\widetilde\partial_{T-}^{\alpha}\eta)z
  \quad\text{strongly in }L^2(0,T;L^2(\Omega)).
\end{equation*}
Along the chosen subsequence, we have
$\mathcal D_\ell^\alpha v_\ell\rightharpoonup w$ and
$v_\ell^+-v_\ell^0\rightharpoonup v-v_0$ in $L^2(0,T;L^2(\Omega))$.
Taking the inner product of \eqref{eq:discrete-caputo-integral-pairing}
with $z$, we pass to the limit by pairing these weakly convergent functions
with the strongly convergent test functions. Since $z$ is arbitrary,
we obtain \eqref{eq:identified-weak-caputo-pairing}.
}
{Since finite sums of functions $\eta(t)z$, with $\eta\in C_c^\infty(0,T)$
and $z\in L^2(\Omega)$, are dense in $L^2(0,T;L^2(\Omega))$,
\eqref{eq:identified-weak-caputo-pairing} determines $w$ uniquely. By the
bound on the discrete derivatives in $L^2(0,T;L^2(\Omega))$, every
subsequence has a further subsequence converging weakly to a limit that
satisfies \eqref{eq:identified-weak-caputo-pairing} by the same argument and
hence equals $w$. Thus $\mathcal D_\ell^\alpha v_\ell\rightharpoonup w$ in
$L^2(0,T;L^2(\Omega))$ along the full sequence.}
\end{proof}

\clearpage
\bibliographystyle{plain}
\bibliography{references}

\end{document}